\documentclass[12pt]{article}
\usepackage{amssymb,latexsym, amsmath, amscd, array, graphicx}
\usepackage{amssymb}

\usepackage{setspace}
\input xy
\xyoption{all}
\makeatletter 

\usepackage{setspace}
\usepackage{enumitem}

\date{}

\newtheorem{theorem}{Theorem}[section]

\newtheorem{proposition}[theorem]{Proposition}

\newtheorem{lemma}[theorem]{Lemma}

\newtheorem{remark}[theorem]{Remark}

\newtheorem{conjecture}[theorem]{Conjecture}

\newcommand{\z}{{\Bbb Z}}
\newcommand{\q}{{\Bbb Q}}
\newcommand{\re}{{\Bbb R}}
\newcommand{\N}{{\Bbb N}}

\newcommand{\map}{{\rm map}}

\newcommand{\invlim}{\raisebox{-1ex}{$\stackrel{\hbox{lim}}{\leftarrow}$}}

\newcommand{\lo}{\rightarrow}
\newcommand{\pq}{\dim_{\z[\frac{1}{q}]} }
\newcommand{\p}{\dim_{\z[\frac{1}{p}]} }

\newcommand{\black}{{\blacksquare}}

\newcommand{\ext}{{\rm ext-dim}}
\begin{document}

\title{\bf Resolving compacta by free $p$-adic actions}

\author{  Michael  Levin\footnote{This research was supported by 
THE ISRAEL SCIENCE FOUNDATION (grant No. 522/14) }}

\maketitle
\begin{abstract}
We say  that a compactum (compact metric space) $Y$ 
is resolvable  by a  $p$-adic action on a compactum $X$ if there is
a continuous  action of the $p$-adic integers $A_p$ on $X$ such that $Y=X/A_p$. 
In this paper we study compacta $Y$ that are resolvable by a  free $p$-adic  action  on a compactum of a lower dimension
and   focus on  compacta $Y$ with $\p Y=1$. This is mainly
 motivated  by   $p$-adic actions on $1$-dimensional compacta, 
the case that turns out to be highly non-trivial. 
More motivation for considering orbit spaces with $\p=1$ comes from 
\\
{\bf   Theorem (A).}  If  $A_p$   acts   on a finite dimensional compactum $X$
 so that $Y=X/A_p$ is infinite dimensional then 
  there exists an invariant compactum $ X'\subset X$  on which 
  the action of $A_p$ is free and whose orbit space $Y'=X'/A_p$ is infinite dimensional 
  with
$\dim_{\z[\frac{1}{p}]} Y' =1$. $\black$

We show 
\\
{\bf Theorem (B).} Let $Y$ be  a finite dimensional compactum with $\p Y =1$.  Then
\\
(i)   $Y$   is resolvable by a free $p$-adic action on 
a  compactum of $\dim \leq \dim Y-1$  if $\dim Y \geq2$ and
\\
(ii) $Y$ is not resolvable by a free $p$-adic action on a compactum of $\dim \leq \dim Y -2$ if $\dim Y\geq 4$.
$\black$

The author  was  initially  inclined to believe that for a  $3$-dimensional compactum  $Y$ with $\p Y=1$
(the case not covered by  (ii) of Theorem (B))  the  additional assumption  $\dim_{\z_p} Y=2$
imposed by Yang's relations would imply that $Y$ is 
 resolvable by a free $p$-adic action on a $1$-dimensional compactum
 and was surprised  to find out that
 \\
 {\bf Theorem (C).}  There is 
 a $3$-dimensional compactum $Y$ with $\p Y =1$  and  $\dim_{\z _p} Y =2$ 
 that cannot be resolved by a free $p$-adic
action on a $1$-dimensional compactum. $\black$

Moreover
\\
{\bf Theorem (D).}  There is 
 an infinite dimensional compactum $Y$ with $\p Y =1$  and  $\dim_{\z } Y =2$ 
 that cannot be resolved by a free $p$-adic
action on a finite dimensional compactum. $\black$

Theorems  (A), (C) and (D) are based on
\\
{\bf Theorem (E). } 
  If  a compactum $Y$  with
  $\dim \geq n+2$  can be resolved by a free $p$-adic action on  an $n$-dimensional compactum 
then  the second Cech cohomology $H^2(Y; \z)$ of $Y$     contains  a subgroup isomorphic to
$\z_{p^\infty}$. Moreover
for every closed $A\subset Y$  with $\dim A \geq n+2$  the image   of this subgroup  in $H^2(A;\z)$
under the homomorphism induced 
by the inclusion is non-trivial. $\black$

\end{abstract}

\begin{section}{Introduction} 

The Hilbert-Smith conjecture asserts that a compact group effectively (and continuously)
acting on a manifold   must be a Lie group.
This assertion is equivalent to the following one:  there is no effective action of  $A_p$ (the group  of  $p$-adic integers) 
on a manifold.  The Hilbert-Smith conjecture is proved for manifolds of $\dim \leq 3$ and open even for
free actions of $A_p$ in $\dim >3$. Yang  \cite{yang}  showed    that if $A_p$ effectively acts on an $n$-manifold $M$ then
either $\dim M/A_p=\infty$ or $\dim M/A_p=n+2$.  This naturally suggests to  examine if the latter dimensional relations
may occur in a more general setting, mainly when $M$  is just
a finite dimensional compactum (=compact  metric space).
One of these relations was  confirmed by 
Raymond and Williams \cite{raymond-williams}  who constructed 
  an action of $A_p$ on an $n$-dimensional compactum  $X$, $n\geq 2$,
with $\dim X/A_p=n+2$. However it remains open    for more than 50 years whether there exists
 a free action of $A_p$ on a finite dimensional compactum $X$ such that
$\dim X/A_p=\dim X +2$ or
$\dim X/A_p=\infty$.

An important collection  of dimensional properties of the orbit spaces under actions of $A_p$ is provided by
  Cohomological Dimension. Recall
the cohomological dimension $\dim_G X$ of a compactum
 $X$ with respect to an abelian group $G$ is 
the least integer  $n$ (or $\infty$ if such $n$ does not exist)  such that the Cech cohomology $H^{n+1}(X,A;G)$
vanishes  for every closed $A\subset  X$. Clearly $\dim_G X  \leq \dim X$ 
 for every group $G$ and,  by the  Alexandrov theorem,
$\dim X=\dim_\z X$ if $X$ is finite dimensional.  A first example of an infinite dimensional compactum  $X$ with 
$\dim_\z  X < \infty$ was constructed by A. Dranishnikov.

 Yang  \cite{yang}  showed  that a free action of $A_p$ on a compactum $X$
imposes the following dimensional relations between $X$ and $Y=X/A_p$:
$\dim_{\z[\frac{1}{p}]} Y=\dim_{\z[\frac{1}{p}]} X$, 
$\dim_{\z_{p^\infty}} X \leq \dim_{\z_{p^\infty}} Y \leq \dim_{\z_{p^\infty}} X +1$, 
 $\dim_{\z_p} X \leq \dim_{\z_p} Y \leq \dim_{\z_p} X +1$  and
$\dim X \leq \dim_\z Y \leq \dim_\z X +2$. Moreover if $\dim_\z Y=\dim_\z X +2$  then
$\dim_{\z_p{^{\infty}}} Y =\dim_{\z_p}  Y =\dim_\z X +1$. We will refer to these dimensional  relations as
the Yang relations.

There is a nice characterization of cohomological  dimension in terms of extensions of maps.
A CW-complex $K$ is said to be an absolute extensor for a  space $X$,   written   \ext$X \leq K$, if 
every map from a closed subset $A $ of $X$ to $K$   continuously extends
over $X$. It turns out that 
 $\dim_G X \leq n$ if and only if  the Eilenberg-MacLane complex $K(G,n)$ is an absolute extensor for $X$.
 Using this characterization  and representing $K(\z[\frac{1}{p}], 1)$ as the infinite telescope of
the $p$-fold covering $S^1 \lo S^1$ of the circle $S^1$ one can easily show 
\begin{proposition}
\label{z[1/p]}
Let $X$ be  a compactum. Then 
$\dim_{\z[\frac{1}{p}]} X \leq 1$  if and only if  for every map $f : A \lo S^1$  from a closed subset 
$A$ of  $X$ to the circle $S^1$ there is a natural number $k$ such that $f$ followed by the $p^k$-fold
covering  map  $S^1 \lo S^1$ of $S^1$  extends over $X$.
\end{proposition}

 Let us
 say that a compactum $Y$ can be  resolved by a free $p$-adic action
on a compactum $X$  if there is a continuous free  action   of the $p$-adic integers $A_p$ on $X$ such that $Y=X/A_p$
 We will  shorten this lengthy expression
to a shorter one: $Y$ is $p$-resolvable by $X$ (keeping in mind that we consider only free actions of $A_p$).

In this paper we study compacta $Y$ that are $p$-resolvable by compacta of a lower dimension
and mainly  focus on compacta with $\p Y=1$ and, in particular on
  compacta   which are   $p$-resolvable by  $1$-dimensional compacta
 (and hence, by Yang's relations, have   $\p  =1$), the case that turns out to be highly non-trivial.
 More motivation for considering orbit spaces $Y$  with $\p Y=1$ comes from
\begin{theorem}
\label{invariant}
 If  $A_p$   acts   on a finite dimensional compactum $X$
 so that $Y=X/A_p$ is infinite dimensional then 
  there exists an invariant compactum $ X'\subset X$  on which 
  the action of $A_p$ is free and whose orbit space $Y'=X'/A_p$ is infinite dimensional 
  with
$\dim_{\z[\frac{1}{p}]} Y' =1$.
\end{theorem}

The most  obvious obstruction to the $p$-resolvability of a compactum $Y$ by  a compactum of 
$\dim < \dim Y$  is $H^1 (Y; \z_p)=0$. Indeed,  $A_p$ is the inverse limit of $\z_{p^n}$ and 
hence every compactum
 that $p$-resolves $Y$
can be represented as the inverse limit of $\z_{p^n}$-coverings of $Y$. 
 Recall that    $H^1(Y; \z_p)$ represents
the $\z_p$-coverings (bundles)  of $Y$. Then  $H^1 (Y; \z_p)=0$  implies that every  $\z_{p^n}$-covering 
over $Y$ is trivial.
 Hence  any compactum
that $p$-resolves $Y$ must contain a copy of every component of $Y$, and therefore cannot be of a lower dimension.

The conditions   $\dim_{\z[\frac{1}{p}]} Y=1$ and $\dim Y >1$ eliminate this obstruction. 
To show that consider
a map from $f : A \lo S^1$ from a closed subset $A $ of $Y$ that does not extend over $Y$. By
Proposition \ref{z[1/p]} there is 
 a $\z_{p^k}$-covering map $\phi :S^1\lo S^1$ such that $f$ followed  by $\phi$ extends
to $g : Y \lo S^1$.  Take a component $Y'\subset Y \times S^1$ of the pull-back space of $g$ and $\phi$
 with the induced maps (projections) $g' : Y' \lo S^1$ and $\phi' : Y' \lo Y$. Notice that $\phi'$ is not $1$-to-$1$ since otherwise
$g'$ would provide an extension of $g$. Thus $\phi'$ is a  $\z_{p^t}$-covering of $Y$ with $0<t\leq k$.
The existence of a non-trivial 
$\z_{p^t}$-covering
of $Y$ 
implies $H^1(Y; \z_p) \neq 0$. 
A similar reasoning also shows that $H^1(Y; \z_p)$ is infinite.

It turns out that  for a finite dimensional compactum $Y$ with $\dim Y >1$ 
 the  condition   $\dim_{\z[\frac{1}{p}]} Y=1$ is already sufficient to
$p$-resolve $Y$ by  a compactum of $\dim < \dim Y$ but not of $\dim < \dim Y -1$.

\begin{theorem}
\label{dim=2}
Let $Y$ be  a finite dimensional  compactum with $\p Y =1$.  Then

(i)   $Y$   is $p$-resolvable by 
a  compactum of $\dim \leq \dim Y -1$  if $\dim Y\geq2$ and

(ii) $Y$ is not $p$-resolvable by  a compactum of $\dim \leq \dim Y -2$ if $\dim Y\geq 4$.
\end{theorem}
The author  was  initially  inclined to believe that  for a $3$-dimensional compactum  $Y$ with $\p Y=1$
(the case not covered by  (ii) of Theorem \ref{dim=3})  the  additional assumption  $\dim_{\z_p} Y=2$
imposed by Yang's relations would imply that $Y$ is  $p$-resolvable by  a $1$-dimensional compactum
 and was surprised  to find out that
  
   \begin{theorem}
\label{dim=3}
There is a $3$-dimensional compactum $Y$ with
$\dim_{\z[\frac{1}{p}]} Y=1$ and $\dim_{\z_p} Y=2$ that is not $p$-resolvable by a $1$-dimensional compactum.
\end{theorem}
Moreover
\begin{theorem}
\label{dim-z=2}
There is an infinite dimensional compactum $Y$ with
 $\dim_{\z[\frac{1}{p}]} Y=1$ and $\dim_\z Y=2$  that is not $p$-resolvable by a finite dimensional compactum.
\end{theorem}
The following property 
 is of crucial importance for  proving  Theorems \ref{invariant}, \ref{dim=3} and \ref{dim-z=2}.
\begin{theorem}
\label{observation} 
 Let  a compactum $Y$  be $p$-resolvable  by  an $n$-dimensional compactum and  $\dim Y \geq  n +2$.
Then  $H^2(Y; \z)$    contains  a subgroup isomorphic to $\z_{p^\infty}$. Moreover
for every closed $A\subset Y$  with $\dim A \geq n+2$ , the image of  this subgroup  in $H^2(A;\z)$ 
under the homomorphism induced 
by the inclusion is non-trivial. 
\end{theorem}
%

 In the proofs    of Theorems \ref{observation},  \ref{dim-z=2} and  \ref{dim=3} 
 we interpret the elements of $H^2(Y; \z)$ as
circle bundles over $Y$  (like  before we  interpreted  the elements of $H^1(Y; \z_p)$ as  $\z_p$-bundles  over $Y$).
A few comments regarding Theorem \ref{observation} are given in  Remark \ref{remark-circle}.

 The author knows  no example of a compactum $Y$   with $\dim Y >2$ and 
 $\p Y=1$  that satisfies the conclusions of
Theorem \ref{observation}.   This together with the fact that the examples in Theorems \ref{dim=3} and
\ref{dim-z=2} are constructed in  a certain generic way motivates

\begin{conjecture}
\label{conjecture}
No  compactum  $Y$ with $\p Y =1$ and $\dim Y \geq n+2$ is $p$-resolvable by an $n$-dimensional compactum.
\end{conjecture}
By  (ii) of  Theorem \ref{dim=3}  the conjecture  holds for  $3 <\dim Y  < \infty$.
Thus the meaningful cases of the conjecture are $\dim Y=3$ and $\dim Y=\infty$.
\end{section}

\begin{section}{Proof of Theorem \ref{dim=2}}

{\bf Kolmogorov-Pontrjagin surfaces.} We will describe $2$-dimensional compacta 
to which we will refer in the sequel  as 
Kolmogorov-Pontjagin surfaces. 

Let $p$ be a prime number and $k \geq 0 $  an integer.
Consider a $2$-simplex $\Delta$, denote by $\Omega(p^k)$ the mapping cylinder of a $p^k$-fold covering map
$\partial \Delta \lo S^1$ and refer to the domain   $\partial \Delta$ and the  range $S^1$ 
of this map  as the bottom and the top of  
$\Omega(p^k)$ respectively. 

Let
  $k_0=0, k_1, k_2,...$ an increasing sequence of natural numbers.
We will construct a Kolmogorov-Pontrjagin surface $Y$ determined by $p$ and  the sequence $k_n$ as the inverse limit 
of $2$-dimensional finite simplicial complexes $\Omega_n$.
Set  $\Omega_0$ to be  a $2$-simplex $\Delta$. 
Assume that $\Omega_n$ is constructed. Take a sufficiently fine triangulation of 
$\Omega_n$ and in every  $2$-simplex $\Delta $ of $\Omega_n$  remove the interior of $\Delta$ and attach to
$\partial \Delta$ the mapping cylinder $\Omega(p^{k_{n+1}})$
 by identifying the bottom of $\Omega(p^{k_{n+1}})$ with
$\partial \Delta$.
  Define $\Omega_{n+1}$ to be  the space that obtained this way from $\Omega_n$   and 
define the bonding map
$\omega_{n+1} : \Omega_{n+1} \lo \Omega_n$ to be a map that sends  each mapping cylinder $\Omega(p^{k_{n+1}})$
to the corresponding simplex   $\Delta$ that identifies the bottom of $\Omega(p^{k_{n+1}})$ with $\partial \Delta$ and
sends the top of $\Omega(p^{k_{n+1}})$ to the barycenter of $\Delta$.
Denote $\Omega=\invlim (\Omega_n, \omega_n)$ and 
call $\Omega$ a Kolmogorov-Pontrjagin surface determined by the prime $p$ and
the sequence $k_n$. We additionally assume that in the construction of a Kolmogorov-Pontrjagin surface $\Omega$
the triangulations of  $\Omega_n$ are so fine that  the diameter of  the images of the simplexes of 
$\Omega_n$  in 
$\Omega_i$, $i < n$ under the map 
$\omega^i_n=\omega_n \circ \dots \circ \omega_{i+1} : \Omega_n\lo \Omega_i$
 is $<1/2^{n-i}$.

\begin{proposition}
\label{maps-to-surfaces}
A compactum $Y$ with $\dim Y \leq m+2$, $m\geq 0$,  and $\p Y=1$ admits an $m$-dimensional
map  into a Kolmogorov-Pontrjagin surface.
\end{proposition}
In the proof of  Proposition \ref{maps-to-surfaces} we will use the following facts.
\begin{proposition} {{\rm \cite{dranishnikov-survey}}}
\label{hurewicz}
 Let 
$f : X \lo Y$ be an  $m$-dimensional map   (= a map whose fibers are of dim $\leq m$).  Then 
 $ \dim X \leq \dim Y +m$ and 
$\dim_G  X \leq  \dim_G Y +m$   for every abelian group $G$. 
\end{proposition}

\begin{lemma}
\label{trees}
Let $T_1$  and $T_2$ be  finite trees,  $X$ a  compactum and $X' $ a  $\sigma$-compact subset of $X$
with $\dim X' \leq 2$.
Then every map $f : X \lo T_1\times T_2$ can be 
arbitrarily closely approximated by a map $f' : X \lo T_1 \times T_2$  such that $f'$ coincides with $f$
on $A=f^{-1}(\partial (T_1 \times T_2))$ and $f'$ is $0$-dimensional on 
$X' \setminus A$ where 
$\partial (T_1 \times T_2)=((\partial T_1 )\times T_2) \cup (T_1\times  \partial T_2)$  and
$\partial T_1$ and $\partial T_2$ stand for the sets of the end points of $T_1$ and $T_2$ respectively.
\end{lemma}
{\bf Proof of Lemma \ref{trees}.} 
Represent $f$ as $f=(f_1, f_2)$  with $f_1 : X \lo T_1$ and $f_2 : X \lo T_2$ being   the coordinate maps.
Fix  metrics on $X$,  $T_1$ and $T_2$, 
consider the induced supremum metric on $Y=T_1 \times T_2$ and let $\epsilon >0$. Since
$T_1$ and $T_2$ are finite trees  there is
$\delta>0$ such for every closed subset $F\subset X$ that does not interesect  $A$
and every  map $f^F_i : F \lo T_i,i=1,2$ 
that  $\delta$-close to $f_i$ on $F$, $f_i^F$  extends  over $X$ to a map that $\epsilon$-close to $f_i$ and
coincides with $f_i$ on $A$. 

Consider a compact subset  $F^\epsilon \subset X'$ such that 
 $F^\epsilon\subset \{ x \in X: d(x, A)\geq \epsilon \}$ and take a finite closed cover $\cal F$ of $F^\epsilon$ 
by subsets of  $F^\epsilon$ of 
diameter$< \epsilon$
such that   the images of the sets in $\cal F$  under both $f_1$ and $f_2$ are of diameter$< \delta$
and $\cal F$
 splits into the union ${\cal F} ={\cal F}_0 \cup {\cal F}_1 \cup {\cal F}_2$ of 
 collections    ${\cal F}_i$ of disjoint sets. 
 
 Denote by $F_i, i=1,2$ 
 the union of  the  sets in ${\cal F}_i$.
 Then   $f_i$ restricted to $F_i$ can be $\delta$-approximated by 
 a map $f^F_i : F_i \lo T_i$ that sends the sets in  ${\cal F}_i$
 to distinct points in $T_i$.
  Now  extend $f^F_i$ to a map  $f^\epsilon_i : X\lo T_i$  that is $\epsilon$-close to $f_i$ and
  coincides with $f_i$ on $A$, and 
   set $f^\epsilon =(f^\epsilon_1, f^\epsilon_2) : X \lo Y=T_1\times T_2$.
   
 Then 
 every fiber of $f^\epsilon$ restricted to $F^\epsilon$  intersects
at most one set of  ${\cal F}_1$ and at most one set  of ${\cal F}_2$ and hence every fiber
of  $f^\epsilon$ restricted to $F^\epsilon$  can be covered 
by finitely many disjoint compact sets of diameter$<4\epsilon$.  
Recall that  $X'$ is $\sigma$-compact, 
 $f^\epsilon$ is $\epsilon$-close to $f$ and coincides with $f$ on $A$, and 
apply the standard Baire category argument 
to the function space of the maps from $X$ to $Y$ that coincide with $f$ on $A$ to get the desired result.
$\black$
\\\\
{\bf Proof of Proposition \ref{maps-to-surfaces}.}
Recall that a Kolmogorov-Pontrjagin surface is the inverse limit of finie simplicial complexes $\Omega_n$ with
the bonding maps $\omega_{n+1} : \Omega_{n+1} \lo \Omega_n$ and 
  $\Omega_0$ is  a $2$-simplex. Take a $\sigma$-compact subset $Y'$ of $Y$ such that $\dim Y'\leq 2$ and
  $\dim Y\setminus Y' \leq m-1$ and 
take any  $m$-dimensional  map $f_0 : Y \lo \Omega_0$.  We will construct by induction on $n$   spaces
$\Omega_n$   and 
$m$-dimensional  maps
$f_n : Y \lo  \Omega_n$. Assume the construction is completed for $n$ and proceed to $n+1$ as follows.
Fix a triangulation of $\Omega_n$.  Recall that  $\dim_{\z[\frac{1}{p}]} Y \leq 1$ and by 
Proposition \ref{z[1/p]} take a sufficiently large $k_n$ such that 
 for every simplex $\Delta$ of $\Omega_n$ 
$f_n$ restricted to  $f^{-1}_n(\partial \Delta)$ and followed by a $p^{k_n}$-covering map
$\partial \Delta \lo S^1$ extends over $f^{-1}_n(\Delta)$. This determines a simlicial complex
$\Omega_{n+1}$ as described in the construction of Kolmogorov-Pontrjagin surfaces and
a natural map $f_{n+1} : Y \lo \Omega_{n+1}$. Clearly taking a sufficiently fine triangulation
of $\Omega_n$ we may assume that $f_n$ and $\omega_{n+1}\circ f_{n+1}$ are as close as we wish.
Note that every point of every $\Omega_i$ has a closed neighborhood homeomorphic to a product of a finite tree
with a closed interval. Then, by Proposition \ref{trees},  we  can replace $f_{n+1}$ by a map which is 
$0$-dimensional on $Y'$. Recall that $\dim Y \setminus Y'\leq m-1$ and we may assume that
$f_{n+1}$ is $m$-dimensional.
Then it is easy to see that  
the whole construction can be carried out  so that the maps 
$f_n$ will determine an $m$-dimensional  map from $Y$ to $\Omega=\invlim \Omega_n$. 
$\black$

\begin{lemma}
\label{covering-omega}
Let $ f :\tilde\Omega(p^k) \lo \Omega(p^k)$ be the $\z_{p^k}$-covering of $\Omega(p^k)$ induced
by the $\z_{p^k}$-covering of the top of  $\Omega(p^k)$. Then $f$ restricted to the preimage of the bottom 
of $\Omega(p^k)$ extends over $\tilde\Omega(p^k) $ as a map to the bottom of $\Omega(p^k)$.

\end{lemma}
{\bf Proof.} Represent $\tilde \Omega(p^k)$ as the union of $p^k$ copies of $S^1\times [0,1]$ being glued 
along  $S^1\times \{ 1 \}$ such that under $f$ the set  $S^1\times \{ 1 \}$   goes to the top of $\Omega(p^k)$
and the sets $S^1 \times \{ 0 \}$ go to the bottom of $\Omega(p^k)$.  Consider a map 
$f_1$ from $ S^1 \times \{ 1 \}$ to the bottom of $\Omega(p^k)$ such that $f_1$ followed by the natural projection
of $\Omega(p^k)$ to its top coincides with $f$ on  $ S^1 \times \{ 1 \}$. Then for each set
$S^1 \times [0,1]$ the maps $f_0=f$ restricted to $S^1\times \{0\}$ and $f_1$ are homotopic as
maps to the bottom of $\Omega(p^k)$ and hence extend over $S^1\times [0,1]$ to  a map 
to the bottom of $\Omega(p^k)$. This defines a map required in the proposition. $\black$

\begin{proposition}
\label{surface-resolvable}
Every  Kolmogorov-Pontrjagin surface is $p$-resolvable by a $1$-dimensional compactum.
\end{proposition}
{\bf Proof.}
Let $\Omega=\invlim \Omega_n$ be a Kolmogorov-Pontrjagin surface determined by a sequence $k_n$.
Consider the first homology $H_1(\Omega_n; \z)$ and let 
 $C_i, i\geq 1$ be  the collection of 
the tops of all the mapping cylinders
added while constructing $\Omega_i, i\leq n$. Then the circles in $C_1\cup \dots \cup C_n$ 
 considered as elements of  $H_1(\Omega_n; \z)$ 
 form a collection of   free generators of  $H_1(\Omega_n; \z)$.

 Take any sequence $s_n$ such that $s_1=k_1$ and $s_{n+1}-k_{n+1} \geq s_n $ for every $n\geq 1$
 and   define  the subgroup $G_n$ of $H_1(\Omega_{n}; \z)$   by
 
 $$  G_n = \{ \sum_{i=1}^n \sum_{\alpha \in C_i} t_\alpha \alpha : 
    \sum_{i=1}^n \sum_{\alpha \in C_i} p^{{s_i}-k_i}  t_\alpha  \ {\rm  is \ divisible \ by} \  p^{s_{n}}\}.$$
 Then 
 $H_1(\Omega_n; \z) /G_n =\z_{p^{s_n}}$
 and for the bonding map $\omega_{n+1} : \Omega_{n+1} \lo \Omega_n$ we have that 
 that the induced homomorphism  $(\omega_{n+1})_* : H_1(\Omega_{n+1}; \z ) \lo H_1(\Omega_n; \z)$
 is onto and sends $G_{n+1}$ into $G_n$. Indeed
 $(\omega_{n+1})_*  (\sum_{i=1}^{n+1}\sum_{\alpha \in C_i} t_\alpha \alpha) =
 \sum_{i=1}^n   \sum_{\alpha \in C_i} t_\alpha \alpha$ and  if
 $\sum_{i=1}^{n+1} \sum_{\alpha \in C_i} p^{s_i-k_i} t_\alpha =
 p^{s_{n+1}-k_{n+1}}(\sum_{\alpha \in C_{n+1}} t_\alpha)+
   \sum_{i=1}^n \sum_{\alpha \in C_i} p^{s_i-k_i}  t_\alpha$ is divisible by $p^{s_{n+1}}$ with
   $s_{n+1}-k_{n+1} \geq s_n$  then we have
   that $\sum_{i=1}^n \sum_{\alpha \in C_i} p^{s_i}  t_\alpha$  is divisible by $s_n$ and therefore
   $(\omega_{n+1})_* (G_{n+1})\subset G_n$.

Let   $h_n : \pi_1(\Omega_n) \lo H_1(\Omega_n; \z)$ be   the Hurewicz homomorphism. 
Then 
$$\tilde G_n=\pi_1(\Omega_n) /h_n^{-1}(G_n)= 
H_1(\Omega_n; \z) /G_n =\z_{p^{s_n}}.$$
 Consider the corresponding $\tilde G_n$-covering
 $f_n : \tilde \Omega_n \lo \Omega_n$ and lift the map $\omega_{n+1} : \Omega_{n+1} \lo \Omega_n$
 to $\tilde \omega_{n+1} : \tilde \Omega_{n+1} \lo \tilde \Omega_n$. The homomorphism 
 $(\omega_{n+1})_*$ induces the natural ephimorpism $\tilde g_{n+1} : { \tilde G_{n+1}} \lo \tilde G_n$
 so that the actions of $\tilde G_{n+1}$ and $\tilde G_n$ 
 agree with   $\tilde g_{n+1}$ and $\tilde \omega_{n+1}$. This defines the compactum
 $X =\invlim (\tilde \Omega_{n}, \tilde \omega_n)$,  the action of
 $A_p =\invlim (\tilde G_n, \tilde g_n)$ on $X$  and the map $f : X \lo \Omega$ determined by the maps $f_n$.
 
Note that if we
 consider   $\alpha \in C_n$ as a circle in $\Omega_n$ then the preimage of $\alpha$ under the map $f_n$ splits
 into $p^{s_n-k_n}$ components (circles)
  and $f_n$ restricted to each component is a $\z_{p^{k_n}}$-covering of  $\alpha$.
 
 Let us show that $\dim X \leq 1$. Consider  the triangulation of $\Omega_n$ used for constructing
 $\Omega_{n+1}$ and a simplex $\Delta$ of this triangulation.
 Let $\Omega(p^{k_{n+1}}) =\omega_{n+1}^{-1}(\Delta)$  and  $\tilde \Omega_{n+1}(p^{k_{n+1}})$ a component
 of $f^{-1}_{n+1} (\Omega (p^{k_{n+1}}))$.  Note that $f_{n+1}$ restricted 
 to $\tilde \Omega(p^{k_{n+1}})$ is a $\z_{p^{k_{n+1}}}$-covering of $\Omega (p^{k_{n+1}})$. Apply 
 Lemma \ref{covering-omega} to extend the map $f_{n+1}$ restricted to 
  the preimage  in $\tilde \Omega(p^{k_{n+1}})$ of the bottom   $\partial\Delta$
 of $\Omega( p^{k_{n+1}})$ to 
  a map $\phi : \tilde \Omega(p^{k_{n+1}})\lo \partial \Delta$.
   Now consider the triangulation of $\tilde \Omega_n$ induced by the triangulation 
  of $\Omega_n$. The map $\phi $ lifts to $\tilde \phi : \tilde \Omega _{n+1}(p^{k_{n+1}}) \lo \tilde \Omega_n$
  so that
   $\tilde \phi (\tilde \Omega _{n+1}(p^{k_{n+1}}) )
   \subset \tilde \omega_{n+1} (\tilde \Omega_{n+1}(p^{k_{n+1}}))$ .
  Doing that for every simplex of $\Omega_n$ we get a map from $\tilde \Omega_{n+1} $
  to the  $1$-skeleton of $\tilde \Omega_n$ and this shows that $\dim X\leq 1$. $\black$
 \\\\
{\bf Proof of Theorem \ref{dim=2}.}

(i)  Let $\dim Y=m+2$. 
 By Proposition \ref{maps-to-surfaces}
there is an $m$-dimensional map $f : Y \lo \Omega$ to a Kolmogorov-Pontrjagin surface $\Omega$. 
By Proposition \ref{surface-resolvable} there is a free  action of  $A_p$ on a $1$-dimensional compactum $\tilde \Omega$
such that $\tilde \Omega /A_p =\Omega$. Let $X$ be the pull-back of $f$ and the projection of $\tilde \Omega$
to $\Omega$.  Then the projection of $X$ to $\tilde \Omega$ is an $m$-dimensional  map and hence,
by Proposition \ref{hurewicz}  $\dim X \leq m+1$ and 
for the pull-back action of $A_p$ on $X$ we have that $X/A_p=Y$. 

(ii)
 Assume that $Y$ is $p$-resolvable by 
a finite dimensional compactum $X$. Then, by Yang's relations $\p X=1$ and  
$\dim_{\z_{p^\infty}} X \geq  \dim_{\z_{p^\infty}} Y -1$. Bokstein inequalities and  $\dim Y \geq 4$
imply that 
$\dim_\q X =\dim _\q Y =1$ and 
$ \dim_{\z_{p^\infty}} Y =\dim Y -1\geq 3 $ and, hence $\dim_{\z_{p^\infty}} X \geq 2$ and 
 $\dim X =\dim_{\z_{p^\infty}} X  +1\geq \dim Y -1$.
 $\black$ 

 \end{section}

\begin{section}{Proof of Theorem \ref{observation}}
A bundle will always mean a locally trivial  principal bundle.  The circle $S^1$ is considered as the group
$S^1=U(1)= \re/\z$. Recall that the classifying space for circle bundles over compacta
is $BS^1=CP^\infty=K(\z,2)$ and therefore every circle bundle  $ Z \lo Y$  over a compactum 
$Y$ is represented by an element $\alpha$ 
of the second Cech cohomology $H^2(Y; \z)$ that being  considered  as
a map $\alpha: Y \lo K(\z,2)$ allows to obtain the bundle $X \lo Y$ as the pull-back of the universal circle bundle
$E \lo K(\z,2)$  and the map $\alpha$. 

Let $Y$ be a compactum and $f : Z \lo Y$ a circle bundle over $Y$. 
  Consider the subgroup $\z_m$ of $S^1$ and in  each fiber of $f$ collapse 
the orbits   under the action of $\z_m$ in $S^1$ to singletons. 
This way we obtain the  circle bundle $f' : Z' \lo Y$
determined by the action of $S^1 =S^1/\z_m$  on $Z'$.
We will refer to $f' : Z' \lo Y$ as the circle bundle  induced by  $f : Z \lo Y$ and $\z_m$.

\begin{lemma} 
\label{circle} 
Let $f : Z \lo Y$ be a circle bundle over a compactum $Y$ represented by $\alpha\in H^2(Y; \z)$ and
let $f' : Z' \lo Y$ be the circle bundle induced by $f$ and $\z_m$. 
Then $f'$ is represented by $\alpha'=m\alpha \in H^2(Y; \z)$.
\end{lemma}
 {\bf Proof.}  Consider the short exact sequence
 $$ 0\lo \z_m \lo S^1 \lo S^1\lo 0. $$
 It defines  a fiber sequence 
 $$ B\z_m \lo BS^1 \lo BS^1 $$
 with the  long exact sequence of the fibration
 
 $$\dots \lo \pi_2 (B\z_m) \lo \pi_2(BS^1) \lo \pi_2(BS^1) \lo \pi_1 (B\z_m) \lo \pi_1(BS^1) \lo \dots$$
Recall that $B\z_m=K(\z_m,1)$ and $BS^1 =K(\z,2)$ and get

$$ 0\lo \z \lo \z \lo \z_m\lo 0. $$
Thus the homomorphism $S^1 \lo S^1$ induces the map $h: BS^1 \lo BS^1$  which acts 
on the second homotopy group of $BS^1$  as the multiplication by $m$. Represent $\alpha$ as
a map $\alpha: Y \lo BS^1=K(\z,2)$. Then $\alpha'$ is represented by $h\circ \alpha$ that translates
in $H^2(Y; \z)$ to $\alpha'=m\alpha$.
$\black$
\\\\
Let $A_p$ freely act on a  finite dimensional
compactum $X$. Consider $A_p$ as a subgroup  the $p$-adic solenoid $\Sigma_p$
and the induced action of $A_p$ on $X \times \Sigma_p$.  Then $\Sigma_p$ naturally acts on 
$X\times_{A_p} \Sigma_p=( X\times \Sigma_p)/A_p$  with $(X\times_{A_p} \Sigma_p) /\Sigma_p=X/A_p$
and   there is a  natural  projection of  $X\times_{A_p} \Sigma_p$ to $S^1=\Sigma_p / A_p$ induced
by the projections  $X \times \Sigma_p \lo \Sigma_p \lo \Sigma_p /A_p$. Note that the fibers of the projection
$X \times_{A_p} \Sigma_p \lo S^1=\Sigma_p / A_p$ are homeomorphic to $X$ and hence, by
Theorem \ref{hurewicz}, $\dim X \times_{A_p} \Sigma_p \leq  \dim X +1$. 

Consider a decreasing sequence
of subgroups $A^k_p$ of $A_p$ 
such that $A^k_p$ is isomorphic  to $A_p$,  $A^0_p =A_p$, $A_p^k /A^{k+1}_p=\z_p$ and
the intersection of all $A^k_p$ contains only $0$.
 Denote  $Y=X/A_p$  and  $Z_n=(X \times_{A_p} \Sigma_p)/A^k_p$.
Note that the circle  $S^1= \Sigma_p / A^k_p$ acts on $Z_k$  with $Y =Z_k /S^1$. This
turns each $Z_k$ into a circle bundle   $f_k :  Z_k \lo Y$ over $Y$  and
the inclusion of $A^{k+1}_p$ into $A^k_p$ defines
 the natural bundle map $g_{k+1} : Z_{k+1} \lo Z_k$
  that  witnesses that 
 the bundle $f_k: Z_k \lo Y$ is induced by 
the bundle $f_{k+1}: Z_{k+1} \lo Y$ and $\z_p$. Moreover, the bundle $f_0: Z_0 \lo Y $ is trivial since
$(X \times_{A_p} \Sigma_p)/A_p =(X/A_p) \times (\Sigma_p /A_p)=Y \times S^1$. 
Also note that $Z=X \times_{A_p} \Sigma_p =\invlim (Z_k, g_k)$  and
the projection $f : Z=X \times_{A_p} \Sigma_p  \lo Y=(X \times_{A_p} \Sigma_p )/\Sigma_p$ coincides with
 the projection of $Z$ to $Z_k$  followed by $f_k$.
 \\
 \\
 {\bf Proof of Theorem \ref{observation}}. Let $A_p$ acts freely on an $n$-dimensional compactum $X$ 
 with $Y=X/A_p$ of $\dim \geq n+2$. Consider the circle bundles  $f_k : Z_k \lo Y$ described above and
 let $\alpha_k \in H^2(Y; \z)$ represent $f_k$.  By Lemma \ref{circle}, $\alpha_k=p\alpha_{k+1}$.
 Recall that $f_0 :  Z_0  \lo Y$ is  a trivial bundle and hence $\alpha_0=0$. Denote
 by $H$ the subgroup of $H^2(Y; \z)$ generated by $\alpha_k, k=0,1, \dots$. We are going to show 
 that $H$ is isomorphic to $\z_{p^\infty}$ and satisfy the conclusions of the theorem.
 
 Let a compactum $A \subset Y$ be of $\dim \geq n+2$.  Aiming at a contradiction  assume
 the image of $H$ under the inclusion of $A$ into $Y$ is trivial in $H^2(A; \z)$. Then every bundle 
 $f_k : Z_k \lo Y$  is trivial over $A$. Let  $\epsilon >0$ be such that $A$ does not admit an open
 cover of  order$\leq n+2$  by sets of diameter$\leq \epsilon$.
   Recall that $\dim Z \leq n+1$. Then, since $Z =\invlim (Z_k,g_k)$, there is a sufficiently
 large $k$ so that $Z_k$ admits an open  cover $\cal U$  of  order $\leq n+2$ so that the images
 of the sets in $\cal U$ under $f_k$ are of diameter$\leq \epsilon$.  Take a section 
 $\phi : A \lo Z_k$ over $A$. Then $\cal U$ restricted to $\phi(A)$ and mapped by $f_k$ to $A$ 
 provides an open cover of $A$ of order$\leq n+2$ by sets of diameter$\leq \epsilon$. Contradiction. 
 
 Now assuming that $A=Y$ we get that  $H$ is non-trivial and, since  $\alpha_{k} =p\alpha_{k+1}$
 and  $\alpha_0=0$,
 we deduce that $H$ is isomorphic to $\z_{p^\infty}$. $\black$
 
 \begin{remark}
 \label{remark-circle}
 \end{remark}

   In general the assumption $\dim Y \geq n+2$ in Theorem \ref{observation}
 cannot be weakened  to $\dim Y \geq n+1$. 
Indeed, it is easy to see that a Kolmogorov-Pontrjagin surface  $\Omega$ has  $H^2(\Omega; \z)=0$  and,
as it was  shown in the previous section,
$\Omega$ is  $p$-resolvable  by a $1$-dimensional compactum.  In this connection we would like to point out
 an iteresting phenomenon  that occurs in the proof of Theorem \ref{observation}
  for $Y=\Omega$: 
all the compacta  $Z_k$ are homeomorphic to $Y \times S^1$
and the bonding maps $g_{k+1} : Z_{k+1} \lo Z_k$
look   the same, 
and still we get that $Z=\invlim Z_k$ is $2$-dimensional despite that each $Z_k$ is $3$-dimensional.
This happens  because we cannot fix  a  trivialization  $Z_k =Y \times S^1$ for each
$Z_k$   so that with respect to these trivializations  each  $g_{k+1}$      will  have  the form
$g_{k+1}(y,s) = (y,ps),  (y,s) \in Y\times S^1$.

However, there is a special case of Theorem \ref{observation} that admits a strengthening. Namely, one can show
  that if   the compactum $Y$  in Theorem \ref{observation}
 is $p$-resolvable by an $n$-dimensional manifold then the assumption $\dim A \geq n+2$
 can be weakened to $\dim A \geq n+1$.

\end{section}

\begin{section}{Proof of Theorem \ref{invariant}}
Let $\mathcal P$ denote the set of all primes. The {\em  Bockstein basis} is the collection of groups
$\sigma=\{ \q, \z_p , \z_{p^\infty}, \z_{(p)} \mid p\in\mathcal P
\}$ where $\z_p =\z/p\z$ is the $p$-cyclic group,
$\z_{p^\infty}={\rm dirlim} \z_{p^k}$  is the $p$-adic circle, and
$\z_{(p)}=\{ m/n \mid n$ is not divisible by $p \}\subset\q$  is the $p$-localization of integers.

By a Moore space $M(G,n)$ for a group  in the Bockstein basis $\sigma$  we always mean the standard Moore space. 
All CW-complexes are assumed to be countable  and  all the spaces are assumed to be separable  metrizable.
A compactum $X$  is said to be hereditarily infinite dimensional if every closed subset of $X$ is either 
$0$-dimensional or infinite dimensional.

\begin{enumerate}[start=1, label={ (\arabic*)}]
\item (Dranishnikov's first extension criterion \cite{dranishnikov-criteria})
\label{dranishnikov-criterion-1}
 Let  $X$ be a compactum  and $K$ a  CW-complex such that \ext$X \leq K$.  
 Then $\dim_{H_n(K)} X \leq n$ for every $n\geq 0$ where
$H_n(K)$ is the reduced  homology of $K$ with the coefficients in $\z$.
\item (Dranishnikov's second extension criterion  \cite{dranishnikov-criteria})
\label{dranishnikov-criterion-2}
 Let  $K$ be a  simply connected  CW-complex, $X$ a finite dimensional  compactum  with
 $\dim_{H_n(K)} X \leq n$ for every $n$ then
\ext$X\leq K$. 

\item(Dranishnikov's splitting theorem  \cite{dranishnikov-splitting, dranishnikov-dydak-splitting})
\label{dranishnikov-splitting}
Let $X$ be a space and $K_1$ and $K_2$ CW-complexes such that 
   \ext$X \leq K_1 * K_2$.   Then $X$ splits into  $X=A_1 \cup A_2$ with \ext$A_1 \leq K_1$, 
\ext$A_2 \leq K_2$ and $A_1$ being $F_\sigma$ in $X$.
\item   \cite{levin-rational}
\label{levin-rational-moore}
 Let  $X$ be a compactum with $\dim_\q X \leq n$.  Then \ext $X \leq M(\q, n)$. 
\item(A factorization theorem that was actually proved in \cite{levin-unstable})
\label{levin-factorization}
 Let $\sigma' $ be a subcollection of  the Bockstein basis $\sigma$ and
 $X$  a compactum (not necessarily finite dimensional)
  such that  \ext $ X \leq  M(G, n_G)$ for every $G \in \sigma'$. 
  Then  every map $f : X \lo K$ to a finite CW-complex $K$ can be arbitrarily closely
approximated by a map $f' : X \lo K$ 
such that  $f'$ factors through a compactum $Z$ with $\dim Z \leq   \dim K$ 
and \ext$ Z \leq M(G, n_G)$ for every $G \in \sigma'$.
\item
\label{ancel}
 It follows from   the  results of  Ancel  \cite{ancel} and  Pol  \cite{ pol} (see also   \cite{levin-inessentiality})  
 on $C$-spaces
that
 every infinite dimensional  compactum  of finite integral cohomological dimension
 contains a hereditarily infinite dimensional compactum. 
 \item \cite{yang}
 \label{yang-action}
 Let $A_p$ act on an $n$-dimensional compactum $X$.
 Then $\dim_\z X/A_p \leq n+3$.  
\end{enumerate}

\begin{proposition}
\label{moore-spaces}
A compactum  $X$  is finite dimensional if and only if there is natural number $n$ such that \ext$X \leq M(\q, n)$ and
\ext$X \leq M(\z_p, n)$ for every prime $p$.
\end{proposition}
{\bf Proof.} If $\dim X=k>0$  then $\dim_G \leq k$ for every group  $G$ and, 
by \ref{dranishnikov-criterion-2}, \ext$X \leq M(G,n)$ for $n=k+1$ and  every $G \in \sigma$.

Now assume that \ext$X \leq M(\q, n)$ and \ext$M(\z_p, n)$ for every prime $p$. Fix $\epsilon>0$  and
 take an $\epsilon$-map (a map with fibers of diameter$< \epsilon$) 
  $f : X \lo K$ to  a finite dimensional cube $K$. 
By  \ref{levin-factorization} $f$ factors through a $(2\epsilon)$-map $g: X\lo Z$ with $Z$ being finite dimensional  
with \ext$Z \leq M(\z_p, n)$ 
for every prime  $p$ and \ext$Z \leq M(\q, n)$. By  \ref{dranishnikov-criterion-1}
$\dim_{\z_p} Z \leq n$ for every prime $p$  and $\dim_\q  Z \leq n$ and hence, by Bockstein inequalities,
$\dim_\z  Z \leq n+1$ and since $Z$ is finite dimensional $\dim Z =\dim_\z   Z \leq n+1$.
Thus for every $\epsilon >0$, $X$ admits a $(2\epsilon)$-map a compactum of
$\dim \leq n+1$ and hence $\dim X \leq n+1$.  $\black$
\\
\\
{\bf Proof of Theorem \ref{invariant}.}
Let $\dim X=n$ and let $A\subset Y$ be the subset of $Y$ corresponding to all the finite orbits.
Then $\dim A \leq n$. Replace  $A$  by a larger a $G_\delta$-subset  of $Y$ with  $\dim A \leq n$.
Since $Y$ is infinite dimensional, we have that $Y\setminus A$ is infinite dimensional, 
and since it  is also   $\sigma$-compact, there is an infinite dimensional compactum in $Y\setminus A$.
Thus replacing $Y$ and $X$  by this compactum and its preimage we may assume that the action of $A_p$ is free.

By \ref{yang-action}  and \ref{ancel}  we may also assume that $Y$ is hereditarily infinite dimensional.
Let us show that $Y$ contains an infinite dimensional compactum of $\p  =1$.  Aiming at a contradiction 
assume that every compactum $Y'$  in $Y$ with $\dim Y' >0$ is of $\p > 1$.  Then
$\pq Y' > 1 $  for every prime $q \neq p$ 
 because otherwise  $H^2( Y'; \z[\frac{1}{q}])=0$ and hence,  by the universal coefficients theorem, 
$H^2(Y'; \z)$ is $q$-torsion and cannot contain a copy of $\z_{p^\infty}$ that violets
Theorem \ref{observation}.

Note that for every prime $q$ we have that the join $M(\z[\frac{1}{q}],1)* M(\z_q, n+3)$ is contractible and hence
\ext$Y \leq M(\z[\frac{1}{q}],1)* M(\z_q, n+3)$. Then, 
by \ref{dranishnikov-splitting},  $Y=Y^0_q \cup Y^1_q $ with \ext$Y^0_q \leq  M(\z[\frac{1}{q}],1)$, 
\ext$Y^1_q \leq M(\z_q, n+3)$ and $Y^0_q$ being $\sigma$-compact.
 By \ref{dranishnikov-criterion-1}, $\pq Y^0_q \leq 1$ and 
 $\dim_{\z_q} Y^1_q \leq n+3$. Since any compactum in $Y$  of positive dimension
 is of $\pq> 1$ we get that $\dim Y^0_q \leq 0$. 
 
 Denote by $Y^0$ the union of $Y^0_q$ for all prime $q$ and let $Y^1=Y\setminus Y^0$.
 Then $\dim Y^0\leq 0$ and \ext$Y^1 \leq M(\z_q, n+3)$ for every prime $q$.   Enlarging $Y^0$ to
 a $0$-dimensional $G_\delta$-subset of $Y$ we get
 that $Y\setminus Y^0\subset Y^1$ is $\sigma$-compact and infinite dimensional and hence 
 contains an infinite dimensional compactum $Y'$  in $Y^1$. 
   By \ref{yang-action},  $\dim_\z Y \leq n+3$ and hence 
 $\dim_\q Y' \leq n+3$.  Then, by \ref{levin-rational-moore}, \ext$Y'\leq M(\q, n+3)$.
 Recall that \ext$Y' \leq M(\z_q, n+3)$ for every prime $q$. Thus, 
  by \ref{moore-spaces}, $Y'$ is finite dimensional, and   we arrive at a contradiction.
 This shows that $Y$ contains an infinite dimensional compactum   with $\p \leq 1$.
 $\black$

\end{section} 
\begin{section}{Auxiliary constructions}
A map between simplicial or CW-complexes is said to be combinatorial if the preimage of every subcomplex of
the range is a subcomplex of the domain.

\label{auxiliary}
\begin{subsection}{Extending partial maps}

\label{partial-maps}
Let $M$ be a finite simplicial complex,
$K$ a connected CW-complex and $f : A \lo K$ a  cellular map  from 
 a subcomplex  $A$ of $M$. We will show how to construct 
a CW-complex $M'$ and  a map $\mu : M' \lo M$ such that $\mu$ restricted to $A'=\mu^{-1}(A)$ and followed by
$f$ extends to a map $f' : M' \lo K$.

Let $A_i= A \cup M^{(i)}$ be the union of  $A$ with the $i$-skeleton $M^{(i)}$ of $M$. 
Extend $f$ over $A_1$ to a cellular map   $f_1: A_1 \lo K$ and  
set $M_1=A_1$ with $\mu_1$ being the inclusion.
Assume that we already constructed a CW-complex $M_i$, cellular maps $\mu_i : M_i \lo A_i$ and
$f_i : M_i \lo K$ so that $f_{i+1}$ extends $f_i$,  $M_i \subset M_{i+1}$ and $\mu_{i+1}|M_i =\mu_i$
for all the relevant indices up to $n$.

The CW-complex $M_{n+1}$ is obtained from $M_n$ by the following procedure.
If $A_n$ already contains $M^{(n+1)}$ set $M_{n+1}=M_n$, $f_{n+1}=f_n$ and $\mu_{n+1}=\mu_n$.
Otherwise,
for every $(n+1)$-simplex $\Delta$ not contained in $A_n$ attach to $M^\Delta_n=\mu^{-1}_n(\partial \Delta)$
the mapping cylinder of $f_n | M^\Delta_n$. The map $f_{n}$ naturally extends over each mapping cylinder
defining $f_{n+1} : M_{n+1} \lo K$. Extend the map  $\mu_{n}$  to $\mu_{n+1}$ by sending each mapping cylinder
to the corresponding simplex $\Delta$  so that the top ($K$-level) of the mapping goes
to the barycenter of $\Delta$ and $\mu_{n+1}$  is linear on the intervals of the mapping cylinder.

Finally for $i=\dim M$ denote $M'=M_i$, $\mu=\mu_i$ and $f' =f_i$.  Note that
$M'$ admits a triangulation for which $\mu$ is combinatorial provided there are  simplicial structures on $M$ and $K$
for which $f$ is simplicial.

\begin{proposition}
\label{isomorphism-circle}
Assume that in the above construction  
$K=S^1$, the subcomplex $A$ contains the  $1$-skeleton of $M$ and  
the map $f: A \lo S^1$ if of degree $p^k$ (that is  $f$
lifts to a $\z_{p^k}$-covering of $S^1$).

 Let $N$ be   a subcomplex   of $M$,   $m=\dim N$ and
    $N'=\mu^{-1}(N)$.
 Consider the homomorphisms 
 $$ (*) \ \  H_m(N' ; \z_{p^t}) \lo H_m(N; \z_{p^t})$$
 
 $$(**) \ \ H_m(N; \z_p) \lo H_m(N; \z_{p^t} ) \lo H_m(M; \z_{p^t} )$$
  
 $$(***)\ \  H_m(N'; \z_p) \lo H_m(N'; \z_{p^{t}} ) \lo H_m(M'; \z_{p^{t}} )$$ 
  induced by the map $\mu$,  the mononorphism $ \z_p \lo \z_{p^t}$  and
   the inclusions of  $N$ into $M$ and 
  $N'$ into $M'$ respectively. 
  Then for $t\leq k$ we have that 
  
  (1) the homomorphism in (*) is an isomorphism and  
  
  (2)
  if the composition of the homomorphisms in (**) is trivial then the composition
   of the homomorphisms  in (***) 
  is trivial as well.

 \end{proposition}
 {\bf Proof.} Note that for $\mu$ is $1$-to-$1$ over the $1$-skeleton of $M$ and for very $2$-simplex
 $\Delta$ of $M$, we have that $\mu$ is either $1$-to-$1$ over $\Delta$ or
 $\mu^{-1}(\Delta)$ is  the mapping cylinder of $f$ restricted $\partial \Delta$
 attached to $\partial \Delta$. Then, since $f$ is of degree $p^k$,  one can easily see 
 that (1) holds  for  $\dim N\leq 2$. 
 
 Now assume  that $\dim N >2$
 and  let $L $ be the union of all the simplexes of $\dim \leq m-3$ of the first barycencetric subdivison of $M$.
 Note that $\dim \mu^{-1}(L)\leq m-2$. Collapse the fibers of $\mu$ over $L\cap N$  and denote by 
 $N^* $ the CW-complex obtained this way  
  from $N'$ and let   $\mu' : N' \lo N^*$ and  $\mu^* : N^* \lo N$ be  the induced natural  maps.

   Let $N_-=N^{(m-1)}$
    and $N_-'=\mu^{-1}(N_-)$.
 Note that $\mu^*$ is $1$-to-$1$  over    $N_-$  and
 hence we can identify $(\mu^*)^{-1}(N'_-)$ with $N_-$. 
Then, since $\dim \mu^{-1}(L)\leq m-2$, we get that
  the map $\mu'$ induces  isomorphisms $H_m(N'; \z_{p^t}) \lo H_m(N^*; \z_{p^t})$
  and $ H_m(N', N'_-; \z_{p^t})\lo H_m(N^*, N_-; \z_{p^t})$.   Moreover,
   it follows from Construction \ref{partial-maps} that,
 since $t \leq k$, for every $m$-simplex $\Delta$ of $N$ (with respect to the origimal triangulation of $M$),
 we have that 
 $H_m({(\mu^*)}^{-1}(\Delta), (\mu^*)^{-1}(\partial \Delta); \z_{p^t})=\z_{p^t}$ 
 and $\mu^*$ induces an isomorphism 
$$H_m({(\mu^*)}^{-1}(\Delta), (\mu^*)^{-1}(\partial \Delta); \z_{p^t})\lo H_m(\Delta, \partial \Delta; \z_{p^t})$$
 (everything here can be easily visualized
 for $m=3$). Then $\mu^*$ also induces an isomorphism 
$ H_m(N^*, N_-; \z_{p^t})\lo H_m(N, N_-; \z_{p^t})$. 
Now  apply 
   the long exact sequences for the pairs $(N, N_-)$ and $(N^*, N_-)$ and the $5$-lemma to show, by induction
   om $m$, 
   that  (1) holds. Note that  we also showed that
\\

\ \ \ ($\dagger$) \ \  $\mu$ induces an isomorphism
  $H_m(N',N_-' ; \z_{p^t}) \lo H_m(N,N_-; \z_{p^t}) $.
  \\

 Let us turn to (2). Clearly we can replace $M$ by its $(m+1)$-skeleton and assume that $\dim M=m+1$.
 Let $M_-$ be the $m$-skeleton of $M$ and
 $M'_- =\mu^{-1}(M_-)$.  By (1) and ($\dagger$)  we get that 
 $\mu$ induces  isomorphisms 
  $H_m(M'_- ;\z_{p^t}) \lo H_m(M_-; \z_{p^t})$ and
 $H_{m+1}(M',M'_-; \z_{p^t}) \lo H_{m+1}(M, M_-; \z_{p^t})$.
  Then applying the long exact sequences of
 the pairs $(M', M'_-)$ and $(M, M_-)$ and the $5$-lemma we get that the homomorphism
 $H_m(M'; \z_{p^t}) \lo H_m(M; \z_{p^t})$ is injective on the image of $H_m(M'_-; \z_{p^t})$
 in  $H_m(M'; \z_{p^t}) $ under the homomorphism induced by the inclusion
 of $M'_-$ into $M'$.  Since $\dim M'_- =\dim N'=m$ we have that the homomorphism
 $H_m(N'; \z_{p^t})\lo H_m(M'_-; \z_{p^t})$ induced by the inclusion of $N'$ into
 $M'_-$  is injective. Combining all this together we get that (2) holds.
 $\black$

 \begin{proposition}
 \label{moore-space}
 Assume that in Construction \ref{partial-maps}    the subcomplex  $A$  contains 
 the $m$-skeleton of $M$ and $K$ is a Moore space $M(\z_{p}, m)$. 
 
  Let $N$ be   a subcomplex   of $M$  with $\dim N \leq m$
 and  $N'=\mu^{-1}(N)$.
 Consider the homomorphisms 
 
  $$ (*) \ \  H_m(N' ; \z_{p^t}) \lo H_m(N; \z_{p^t})$$
 $$(**) \ \ H_m(N; \z_p) \lo H_m(N; \z_{p^t} ) \lo H_m(M; \z_{p^t} )$$
  
 $$(***)\ \  H_m(N'; \z_p) \lo H_m(N'; \z_{p^{t+1}} ) \lo H_m(M'; \z_{p^{t+1}} )$$ 
  induced by the mononorphisms $ \z_p \lo \z_{p^t}$ and $\z_p \lo \z_{p^{t+1}}$ and
   the inclusions of  $N$ into $M$ and 
  $N'$ into $M'$ respectively.  Then for every $t>0$ we have that 
  
  (1) the homomorphism in (*) is an isomorphism and  
  
  (2)
  if the composition of the homomorphisms in (**) is trivial then the composition
   of the homomorphisms  in (***) 
  is trivial as well.

 \end{proposition}
{\bf Proof.} Note that $\mu$ is $1$-to-$1$ over the $m$-skeleton of $M$
and hence we can identify the $m$-skeleton $M^{(m)}$  of $M$ 
with $\mu^{-1}(M^{(m)})$  and assume that $N=N'$. In particular it implies (1).
We also can replace $M$ by the $(m+1)$-skeleton of $M$ and assume that $\dim M \leq m+1$.
Consider  $\alpha \in H_m(N;  \z_p)$  and let $\beta\in H_m(N; \z_{p^t})$ 
be the image of $\alpha$ and $\gamma\in H_m(N; \z_{p^{t+1}})$  the image
of $\beta$ under  the homomorphisms 
$ H_m(N; \z_p) \lo H_m(N; \z_{p^t} )$  and
$ H_m(N; \z_{p^t} ) \lo H_m(N; \z_{p^{t+1}} )$
 induced by the monomorphisms
$\z_p \lo \z_{p^{t}} $  and $ \z_{p^{t}} \lo \z_{p^{t+1}}$ respectively.

Consider  $\beta$ as an 
 $m$-cycle in the chain complex 
 $C(N; \z_{p^{t}})$. Then $\gamma$ is  represented by the cycle $p\beta$ in the chain complex
 $C(N; \z_{p^{t+1}})$
 Since the  composition (**) is trivial,   the cycle $\beta$
 is homologous to $0$ in $C(M; \z_{p^t})$ and hence $\beta=\partial \theta$ where
 $\theta=c_1 \Delta_i + \dots +c_k\Delta_k$, $c_i \in \z_{p^t}$ and $\Delta_i$ are $(m+1)$-simplexes in $M$.
 Thus $\gamma =p\partial \theta=c_1(p\partial \Delta_1)+ \dots +c_k (p\partial \Delta_k)$
 in  $C(N;\z_{p^{t+1}})$. Recall that 
 for every $(m+1)$-simplex $\Delta$, $\mu^{-1}(\Delta)$ is either contractible or homotopy equivalent
 to the mapping cylinder of $f|\partial \Delta \lo M(\z_{p}, m)$ and hence $p\partial \Delta$ is homologous
 to $0$ in $C(M';\z_{p^{t+1}}) $. Thus we get that  $\gamma$ is homologous to $0$ in $C(M';\z_{p^{t+1}}) $
 and hence
 the composition (***)   is trivial. $\black$

\end{subsection}

\begin{subsection}{Flexible circle bundles}
\label{flexible-bundles}
Let $M$ be a finite simplicial complex and $g : L \lo M$ a circle bundle. We will say that $g$ is $p$-flexible
if for every $k >0$ there is a section $ f : M^{(1)} \lo L$ over the $1$-skeleton of $M$ such that
for a every $2$-simplex $\Delta$ of $M$ and a trivialization $g^{-1}(\Delta)=\Delta\times S^1$ we have  that
the map $f$ restricted to $\partial \Delta$ and followed by the projection $\Delta \times S^1 \lo S^1$
lifts to the $ \z_{p^k}$-covering of $S^1$. We will refer  to $p^k$ as a degree of $f$.
Clearly this definition does not depend on the trivializations
over $2$-simplexes of $M$.

We will describe two versions  of killing non-trivial $p$-flexible bundles that will be used in different contexts.

\begin{subsubsection}{Killing non-trivial flexible bundles: a version for Theorem \ref{dim=2}}
\label{flexible-1}
Let $ g : L \lo M$ be a $p$-flexible circle bundle over a finite simplicial complex $M$. 
Let $f : M^{(1)} \lo L$ be a section of degree $p^{k}$ over the $1$-skeleton of $M$. We will construct
a CW-complex $M'$ and a map $\mu : M' \lo M$ so that  the pull-back bundle $ g': L' \lo M'$ of $f$ via $\mu$
admits a section $f' : M' \lo L'$. 

Set $M_1=M^{(1)}$, $\mu_1: M_1 \lo M$ to be  the inclusion and $f_1 =f: M_1 \lo L$.
We will construct by induction finite CW-complexes $M_i$, maps $\mu_i : M_i \lo M^{(i)}$
 and $f_i : M_i \lo L$ so that $g \circ f_i =\mu_i$.

 Assume that the construction is completed for $i \leq n$
 and proceed to $n+1$  as follows.  Consider $(n+1)$-simplex $\Delta$ of $M$ consider a trivialization
 $g^{-1}(\Delta)=\Delta \times S^1$ and let    $f_\Delta : \mu_n^{-1}(\partial \Delta) \lo S^1$ be
 the map $f_n$  restricted to $\mu_n^{-1}(\partial \Delta)$  and 
  followed by the projection  $g^{-1}(\Delta)=\Delta \times S^1 \lo S^1$.
 
 Consider the mapping cylinder $M_\Delta$ of $f_\Delta$ 
 and let $\mu_\Delta : M_\Delta \lo \Delta$ be the map sending the top of $M_\Delta$ to
 the barycenter of $\Delta$ and linearly extending $\mu_n$ restricted to $\mu_n^{-1}(\partial \Delta)$.
 Then the map $f_n$ restricted  $\mu_n^{-1}(\partial \Delta)$ naturally  extends over   $M_\Delta$ to
  a map $f_\Delta : M_\Delta \lo g^{-1}(\Delta)$   such that  $\mu_\Delta =g \circ  f_\Delta $.
Now  attach  $M_\Delta$ to  $\mu_i^{-1}(\partial \Delta)$ 
 and doing that for
 every $(n+1)$-simplex  $\Delta$  of $M$ we obtain the CW-complex $M_{n+1}$ and the maps
 $\mu_{n+1} : M_{n+1} \lo M^{(n+1)}$ and $f_{n+1} : M_{n+1} \lo L$ induced by the maps $\mu_\Delta$ and
 $f_\Delta$ respectively. Clearly  $g \circ f_{n+1}=\mu_{n+1}$.

Finally for $i=\dim M$ set  $M'=M_i, \mu' =\mu_i$, $ L'=$ the pull-back of $L$ under $\mu'$,  $g$ and $g' : L'\lo M'$
the pull-back of $g$. Then $f_i$ induces a section $f' : M' \lo L'$. Note that
$M'$ admits a triangulation for which $\mu$ is combinatorial.

\begin{proposition}
\label{bundle}  The conclusions of Proposition \ref{isomorphism-circle} hold for the  construction  above.
 \end{proposition}
The proof of Proposition \ref{isomorphism-circle}  applies to prove  Proposition \ref{bundle}  as well.
\end{subsubsection}

\begin{subsubsection}{Killing non-trivial flexible bundles: a version for Theorem \ref{dim-z=2}}
\label{flexible-2}
Let $ g : L \lo M$ be a $p$-flexible circle bundle over a finite simplicial complex $M$. Let $m=\dim M$ and $k=tm$, and
let $f : M^{(1)} \lo L$ be a section of degree $p^{km}$ over the $1$-skeleton of $M$. We will construct
a CW-complex $M'$ and a map $\mu : M' \lo M$ so that  the pull-back bundle $ g': L' \lo M'$ of $f$ via $\mu$
admits a section $f' : M' \lo L'$. 

Set $M_1=M^{(1)}$, $\mu_1: M_1 \lo M$ to be  the inclusion and $f_1 =f: M_1 \lo L$.
We will construct by induction finite CW-complexes $M_i$, maps $\mu_i : M_i \lo M^{(i)}$
 and $f_i : M_i \lo L$ so that $g \circ f_i =\mu_i$ and 
 $g_i$ is a map of degree $p^{(m-i)k}$. By this we mean that   for every $(i+1)$-simplex $\Delta$ of $M$
 and a trivialization $f^{-1}(\Delta)=\Delta \times S^1$ of $g$ over $\Delta$ the map $f_i$ restricted to
 $\mu_i^{-1}(\partial \Delta)$ and followed by the projection $\Delta \times S^1 \lo S^1$  lifts 
 to the $\z_{p^{(m-i)k}}$-covering of $S^1$.
 
 Assume that the construction is completed for $i \leq n$
 and proceed to $n+1$  as follows.  Consider $(n+1)$-simplex $\Delta$ of $M$ consider a trivialization
 $g^{-1}(\Delta)=\Delta \times S^1$ and a map  $f_\Delta : \mu_n^{-1}(\partial \Delta) \lo S^1$
 witnessing that $f_n$ restricted to  $\mu_i^{-1}(\partial \Delta)$ is a map of degree $p^{(m-n)k}$ and
 consider the  finite telescope   $T$ of  $m$ maps $S^1 \lo S^1$  each of them is    the $\z_{p^t}$-covering of $S^1$.  
  the mapping cylinder $M_\Delta$ of $f_\Delta$ followed by the embedding  of $S^1$ in $T$ as the first circle  of $T$.
 Let $\mu_\Delta : M_\Delta \lo \Delta$ be the map sending the top ($T$-level) of $M_\Delta$ to
 the barycenter of $\Delta$ and linearly extending $\mu_n$ restricted to $\mu_n^{-1}(\partial \Delta)$.
 Then the map $f_n$ restricted  $\mu_n^{-1}(\partial \Delta)$ naturally  extends over   $M_\Delta$ to
  a map $f_\Delta : M_\Delta \lo g^{-1}(\Delta)$  of degree $p^{(m-n-1)k}$ such that
   $\mu_\Delta = g \circ f_\Delta$.
Now  attach  $M_\Delta$ to  $\mu_i^{-1}(\partial \Delta)$ 
 and doing that for
 every $(n+1)$-simplex  $\Delta$  of $M$ we obtain the CW-complex $M_{n+1}$ and the maps
 $\mu_{n+1} : M_{n+1} \lo M^{(n+1)}$ and $f_{n+1} : M_{n+1} \lo L$ induced by the maps $\mu_\Delta$ and
 $f_\Delta$ respectively. Clearly  $\mu_{n+1}=g \circ f_{n+1}$. 
 
 The only thing that we need 
 to check is that  that $f_{n+1}$ is of degree  $p^{(m-n-1)k}$ on $\mu_{n+1}^{-1}(\Delta)$ for
 every $(n+2)$-simplex $\Delta $ of $M$.

  The only thing that we need 
 to check is that  that $f_{n+1}$  restricted  to  $\mu_{n+1}^{-1}(\partial \Delta)$ is of degree  $p^{(m-n-1)k}$  for
 every $(n+2)$-simplex $\Delta $ of $M$.  Consider a trivialization  
 $g^{-1}(\Delta)=\Delta \times S^1 $ of $g^{-1}(\Delta)$ and denote by $\Delta_0, \dots, \Delta_{n+1}$
 the $(n+1)$-simplexes  contained in $\Delta$.  Note that for 
 $\mu_{n+1}^{-1}((\Delta_0 \cup \dots \cup \Delta_i) \cap \Delta_{i+1})$ is connected for every
 $0 \leq i \leq n$. Then the cycles of $H_1(\mu_{n+1}^{-1}(\Delta_i))$, $0\leq i \leq n+1$ generate
 $H_1(\mu_{n+1}^{-1}(\partial \Delta))$. Since $f_{n+1}$ restricted to each $\Delta_i$ and followed by
 the projection $\alpha: \Delta \times S^1 \lo S^1$ is of degree $p^k$ we get that
the homomorphism induced by $\alpha \circ f_{n+1}$ sends
$H_1(\mu_{n+1}^{-1}(\partial \Delta))$ into the elements of   $H_1(S^1)=\z$
divisible by $p^{(m-n-1)k}$.  This means
that   $f_{n+1}$  restricted  to  $\mu_{n+1}^{-1}(\partial \Delta)$ is of degree  $p^{(m-n-1)k}$.

Finally set $m=\dim M$, $M'=M_m, \mu' =\mu_m$, $ L'=$ the pull-back of $L$ under $\mu'$,  $g$ and $g' : L'\lo M'$
the pull-back of $g$. Then $f_m$ induces a section $f' : M' \lo L'$. Note that
$M'$ admits a triangulation for which $\mu$ is combinatorial.

Let $T^j, 1\leq j \leq m$ be the subtelescope of $T$  consisting of the first $j$ maps  $S^1 \lo S^1$.
 Thus $T^m=T$ and 
in the construction above we can consider  the  subcomplexes  $M^j \subset M'=M^m$ obtained from $M'$ by 
leaving only the subtelescope $T_j$ from $T$ in all the mapping cylinders $M_\Delta$.

\begin{proposition} 
\label{flexible-2-more}
Let $N$ be  a connected CW-complex whose homotopy groups are finite $p$-torsion groups.
Consider Construction \ref{flexible-2} with 
$t$ being such  such that $p^t\pi_n(N)=0$  for every $1\leq n \leq m$ and
assume  that  $f_N : A_N \lo N$ is a map from a subcomplex $A_N$ of $M$  such that $f_N$ does not extend over $M$.
Then $\mu$ restricted to $\mu^{-1}(A_N)$ and followed by $f_N$ does not extend over $M'$.
\end{proposition}
 {\bf Proof.} Let  $M_*$ by the space obtained from $M'$ by collapsing  the fibers of $\mu$ over $N$ to singletons,
$\mu_* : M_* \lo M$ the  map induced by $\mu$ and $M^j_* =\mu_*(M^j), 1\leq j \leq m$
Recall that  $M^m=M'$ and hence $M_*^m=M_*$.  
Thus we can consider $N$ as a subcomplex of $M_*$ and 
aiming at a contradiction assume that  $f_N$  extends
to $f_* : M_*  \lo N$.  Recall that in Construction \ref{flexible-2} we denote by $T$ the telescope of
$m$ copies of $\z_{p^t}$-coverings  $S^1 \lo S^1$ and by $T^j$, $ 1\leq j \leq m$  the subtelescope of the first $j$ 
maps of $T$ ($T^m=T$).

Consider the first barycentric subdivision $\beta M$ of the triangulation of $M$.
Note that for ever $0$-dimensional simplex (vertex) $\Delta$ of $\beta M$,
$(\mu_*)^{-1}(\Delta)$  is either a singleton or homeomorphic to $T=T^m$. Then
 $(\mu_*)^{-1}(\Delta) \cap M^{m-1}_*$ is either a singleton
or homeomorphic  to the subtelescope $T^{m-1}$ of  $T$ and hence 
$f_*$ is null-homotopic on $(\mu_*)^{-1}(\Delta) \cap M^{m-1}_*$.
 Thus we can replace $M_*^{m-1}$ by the space  $M^{m-1}_0$  obtained from $M_*^{m-1}$ by collapsing
$(\mu^*)^{-1}(\Delta) \cap M^{m-1}_*$  to singletons for every
$0$-simplex $\Delta$ of $\beta M$,  and $\mu_*$ and $f_*$  by the induced maps $\mu^{m-1}_0$ and $f^{m-1}_0$
to $M$ and
$N$ respectively   and assume that
$\mu_0^{m-1}$ is $1$-to-$1$ over the $0$-simplexes of $\beta M$.

Now consider a $1$-simplex $\Delta$ of $\beta M$. Then $(\mu_0^{m-1})^{-1}(\Delta)$ is either contractible 
or homotopy equivalent to $\Sigma T^{m-1}$.  Then
 $(\mu_0^{m-1})^{-1}(\Delta) \cap M^{m-2}_0$ is either a singleton
or homeomorphic  to the subtelescope $\Sigma T^{m-2}$ of  $\Sigma T^{m-1}$
and hence  $f^{m-1}_0$ is null-homotopic on
$(\mu^{m-1}_0)^{-1}(\Delta) \cap  M^{m-2}_0$. 
Thus $f^{m-1}_0$ restricted to $M^{m-2}_0$ factors up to homotopy  through the space obtained from
$M^{m-2}_0$ by collapsing the fibers of $\mu^{m-1}_0|M^{m-2}_0$ over the simplex $\Delta$. 
Doing that consecutively 
for all the  $1$-simplexes of $\beta M$  we obtain the space $M^{m-2}_1$ and 
 the maps $\mu^{m-2}_1 : M^{m-1}_1 \lo M$
and $f^{m-2}_1 : M^{m-2}_1 \lo N$  induced by $\mu^{m-1}_0$ and $f^{m-1}_0$ respectively such that 
$\mu^{m-2}_1$ is $1$-to-$1$ over the $1$-simplexes of $\beta M$.

Procced by induction  and construct for every
$i\leq m-1$ the space $M^{m-i-1}_i$ and the maps $\mu^{m-i-1}_i : M^{m-i-1}_i  \lo M$
and $f^{m-i-1}_i : M^{m-i-1}_i \lo N$
and finally get for $i=m-1=\dim M  -1$ that $M^{0}_{m-1}=M$ and
$f^0_{m-1}$ extends $f_N$ that  contradicts the assumptions of the proposition. $\black$

\end{subsubsection}

\end{subsection}

\begin{subsection}{Constructing  inverse limits}
\label{inverse-limit}

We will describe  how to construct
 an inverse limit $Y=\invlim (M_i, \mu_i)$  of finite simplicial complexes $M_i$   performing
 countably many 
 procedures  on certain objects determined by each $M_i$.
   We will mainly deal with  the objects     described in Section \ref{auxiliary} 
 like partial maps and and circle bundles and the procedures like extending partial partial maps  and killing  some
 non-trivial circle  bundles.
  We assume  that  any object  on  $M_j$ can be transferred to any $M_i$ with $i\geq j$ 
  via  the map $\mu_i^j =\mu_i \circ \dots \circ \mu_{j+1}: M_i \lo M_j$ ($\mu^i_i$ is  the identity  map of $M_i$).
By transferring   a partial map on $M_j$ to $M_i$
 we just mean that a partial map $f : A \lo K$ from a closed subset  $A$ of $M_j$ to 
  to CW-complex $K$ moves
 to
 the partial map $f \circ \mu_i^j : (\mu_i^j)^{-1} (A) \lo K$  on $M_i$. 
  And by transferring a circle bundle over $M_j$
 to $M_i$ we mean the pull-back of the bundle  to $M_i$ via the map $\mu^j_i$. 
 We also assume for each  finite simplicial complex 
  we  can fix   countably many  objects on which we want to perform   appropriate procedures.

The inverse system $(Y_i, \mu_i)$  is constructed as follows.
Consider a bijection $\beta : \N \times \N  \lo \N$ such that $\beta(j, n)\geq j$.
Assume that $Y_j$ and $\mu_j$ are already constructed for $j \leq i$. Moreover
  for each $j\leq i$ we also fixed
countably many objects  ${\cal Q}_j$     for $Y_j$ that we need to take care of   and the objects
in ${\cal Q}_j$ are indexed by natural numbers.
Proceed to $i+1$ as follows. Let $(j,n)=i$. Transfer the object indexed by $n$ in ${\cal Q}_j $ to $M_i$, 
construct  $M_{i+1}$ and $\mu_{i+1}$  in order to perform  the  procedure  appropriate 
 for this object. And, finally,
pick out countably  many  objects  ${\cal Q}_{i+1}$   for $Y_{i+1}$ needed to be taken care of  and index
the objects in ${\cal Q}_{i+1}$  by the natural numbers.

Thus we get in the inverse limit determining  $Y$  we took care of  all the objects  picked out for each $M_i$.

\end{subsection}
\end{section}
\begin{section}{Proof of Theorem \ref{dim=3}} 
Theorem \ref{dim=3} follows form Theorem \ref{observation} and the case $n=1$ of the following
proposition.
\begin{proposition}
For every $n \geq 1$ there is a compactum $Y$ with $\p Y=1$, $\dim_{\z_p}  Y =n+1$ and
$\dim Y =n+2$ such that $H^2 (Y; \z)$ does not contain a subgroup isomorphic to $\z_{p^\infty}$.
\end{proposition}
{\bf Proof.}
We will construct   $Y$ as the inverse limit of $(n+2)$-dimensional  finite  simplicial complexes $M_i$ and
 combinatorial
 bonding maps $\mu_{i+1}: M_{i+1} \lo M_i$. In order to show that 
 $Y$  has the required properties
 we consider
 for each  $i$   a subcomplex $A_i$ of
 $M_i$ such that  $A_{i+1}=\mu^{-1}_{i+1}(A_i)$ and 
  two natural numbers $k_i$ and $t_i$ 
 such that $k_{i+1} \geq k_i$, $t_{i+1}\geq t_i$ and $ k_i \geq t_i$.
Set  $M_0$ to be an  $(n+2)$-ball, $A_0=S^{n+1}$ the boundary of this  ball,
 and  $k_0=t_0=1$.
 We require that

\begin{enumerate}[start=1, label={ (\arabic*)}]

\item $\dim A_i=n+1$, $H_{n+1} (A_i; \z_{p}) =\z_p$ and the homomorphism 
$H_{n+1} (A_{i+1}; \z_p) \lo H_{n+1}(A_i; \z_p)$ induced by $\mu_{i+1}$ is an isomorphism;

\item 
Let   $H_{n+1}(A_{i}; \z_{p})  \lo H_{n+1} (A_i; \z_{p^{t_i}})$ and 
 $H_{n+1} (A_i; \z_{p^{t_i}})\lo H_{n+1} (M_i; \z_{p^{t_i}})$ be the homomorphisms 
  induced by the monomorphism
 $\z_p \lo \z_{p^t}$ and the inclusion of $A_i$ into $M_i$ respectively. Then
the following composition is trivial
 $$H_{n+1}(A_{i}; \z_{p})  \lo H_{n+1} (A_i; \z_{p^{t_i}})\lo H_{n+1} (M_i; \z_{p^t}).$$

\end{enumerate}
Clearly  the relevant parts of (1)  and  (2)  hold for $i=0$.
Assuming that  the construction is completed for $i$  and we proceed to $i+1$ performing one the following 
procedures.\\

 {\bf Procedure I (taking care of non-flexible  bundles).} Let $g : L \lo M_i$ be a circle bundle which is not
   $p$-flexible. Then take  any  natural number $k$ such that  there is no section over $M_i^{(1)}$ of degree $p^k$.
   Set $M_{i+1}=M_i$, $\mu_{i+1}$=the identity map, $k_{i+1}= \max\{k, k_i\}$ and $t_{i+1}=t_i$.\\
   
   {\bf Procedure II (taking care of flexible bundles).}
   Let $g : L \lo M_i$ be a flexible bundle. Apply \ref{flexible-1} with $M=M_i$ and $k=k_i$ 
   construct $M_{i+1}=M' $ and $\mu_{i+1}=\mu$.  Set $k_{i+1}=k_i$ and $t_{i+1}=t_i$.\\
   
{\bf Procedure III (extending partial maps to a circle).} 
Let $A$ be a subcomplex of $M_i$ and $f : A \lo S^1$ a map.
Extend $f$ over  the $1$-skeleton $M^{(1)}$ of $M$  and replace $f$ by  a cellular map
homotopic  to $f$ followed by a map  $S^1 \lo S^1$
of degree $k_i$. Thus we assume that  $A$ contains the $1$-skeleton of $M_i$ and
  $ f : A \lo S^1$ is  a cellur map of degree $p^{k_i}$
   (recall that   ``of degree $p^{k_i}$'' that
    $f$ lifts to a $\z_{p^{k_i}}$-covering of $S^1$). Apply \ref{partial-maps} with $M=M_i$ and
   $K=S^1$  
  to  construct $M_{i+1}=M' $ and $\mu_{i+1}=\mu$.  Set $k_{i+1}=k_i$ and $t_{i+1}=t_i$.\\
   
  {\bf  Procedure  IV (extending partial maps to $M(\z_p, n+1)$).} 
  Let $A$ be a subcomplex of $M_i$
  and $f: A \lo M(\z_p, n+1)$ a map. Extend $f$ over the $(n+1)$-skeleton on $M_i$ and
  replace $f$ by a homotopic cellular map.  Thus we assume that $A$ 
   contains the $(n+1)$-skeleton $M^{(n+1)}$ of $M_i$ and 
   $f : A \lo M(\z_p, n+1)$ is  a cellular map.
  Apply \ref{partial-maps} with $M=M_i$ and $K=M(\z_p, n+1)$ to
   construct $M_{i+1}=M' $ and $\mu_{i+1}=\mu$.  Set $k_{i+1}=k_i+1 $  and $t_{i+1}=t_i +1$.
   \\
   \\
Let us check that after preforming the above procedures the conditions (1)  and   (2) hold.
It is obvious for Procedure I.  Procedures II and III preserve (1) and (2) by Propositions \ref{isomorphism-circle}
and \ref{bundle} because $t_{i+1}=t_i \leq k_i$. Procedure IV preserves (1) and (2) by Proposition \ref{moore-space}.

Let us show that $Y=\invlim M_i$ is $(n+2)$-dimensional. Clearly $\dim Y \leq n+2$.
Consider the map $\mu_i|A_i : A_i \lo A_0=S^{n+1}$. 
 Since  $\dim A_{i}\leq n+1$,  the long exact sequence generated by
$$0 \lo \z_{p} \lo \z_{p^{t_i}} \lo \z_{p^{t _i-1}}$$ implies
 that $H_{n+1}(A_{i}; \z_{p}) \lo H_{n+1}(A_{i}; \z_{p^t})$ is injective for every $i$ and
 hence, by (1) and (2),  $\mu_i$ induces a non-trivial homomorphism  
 
 $$ \ker( H_{n+1}(A_i; \z_{p^{t_i}}) \lo H_{n+1}(M_i, \z_{p^{t_i}}) )
 \lo H_{n+1}(A_0; \z_{p^{t_i}}).$$
 This implies that $\mu_i$ restricted to $A_i$ does not extend over $M_i$ as a map to $A_0=S^{n+1}$
 and hence $\dim Y=n+2$.
 
 Assume that $Y$ is constructed as described in \ref{inverse-limit}. 
 The objects that we pick out for each $M_i$ are
 all the  (non-isomorphic) circle bundles  over $M_i$ and 
 countably many maps from subcomplexes of $M_i$ to $S^1$ and $M(\z_p, n+1)$ representing
 all possible maps up to homotopy.  The procedures that we perform are Procedures I-IV.
 
 Assuming that on each $M_i$ we fix a sufficiently finite triangulation, we  get that:
 
 $\bullet$  Procedure III together with 
 Proposition \ref{z[1/p]}  implies   that $\p Y \leq 1$ and
 
 $\bullet$ Procedure IV implies  that
 \ext$ Y \leq M(\z_p, n+1)$ and hence, by Dranishnikov's extension criterion, $\dim_{\z_p} Y \leq n+1$.
 
 Then $\dim Y=n+2$ and Bockstein inequalities
 imply that $\p Y =1$ and $\dim_{\z_p} Y={n+1}$.
 
  Let us show that $H^2(Y; \z)$ does not contain a subgroup isomorphic to $\z_{p^\infty}$. Since
  $\p Y =1$ we have that  $H^2(Y; \z)\otimes \z[\frac{1}{p}]=0$ and hence   $H^2(Y; \z)$ is $p$-torsion.
  Thus we need to show  $H^2(Y; \z)$ does not contain a non-trivial
   element $\alpha$ that is divisible by $p^k$ for every $k$.
  Aiming at a contradiction assume that such $\alpha$ does exist. Consider the circle bundle
  $g_Y : Z \lo Y$ corresponding to $\alpha$. Then Construction \ref{inverse-limit}
  guarantees that at a certain step $i$
of the construction we will consider a circle bundle $g : L\lo M_i$ such that
$g_Y$ is the pull-back of $g$ under the projection $\mu^i: Y \lo Y_i$.

 If $g$ is flexible then we apply Procedure II and,
by Construction \ref{flexible-1}, 
the pull-back bundle $g_{i+1} : L_{i+1} \lo M_{i+1}$ of $g$ under the map $\mu_{i+1}$ admits
 a section and hence both $g_{i+1}$  and  $g_Y$ are trivial. Thus we arrive at a  contradiction with $\alpha\neq 0$.
 
 If $g: L \lo M_i$ is non flexible then we apply Procedure I. Denote by $g_j : L_j \lo M_j, j\geq i$, 
 the pull-back bundle of $g$   to $M_j$ under 
 $\mu^i_j =\mu_j \circ \dots \circ \mu_{i+1}: M_j \lo M_i$ with $\mu^i_i$ being the identity map
 and   $g_i=g$ and $L_i=L$. Denote by $G^i_j : L_j \lo L_i$ the induced map.
 Recall that  $\alpha$ is divisible by $p^k$  for every  $k$. Then for a sufficiently large $j$
  there is a circle bundle $g^*_j : L^*_j \lo M_j$ such that $g_j : L_j \lo M_j$ is induced  by 
  $g^*_j$ and $\z_{p^{k_{i+1}}}$,  see Lemma \ref{circle}.  Let 
  $G_j^* : L^*_j \lo L_j$ be a bundle map induced by the homomorphism $S^1 \lo S^1 / \z_{p^{k_{i+1}}}$,  
Take a section $s^*_j : M^{(1)}_j \lo L^*_j$.

  Consider a $2$-simplex $\Omega_i$ of $M_i$ and let $g^{-1}_i(\Omega_i)=\Omega_i \times S^1$ be
  a trivialization of $L_i$ over $\Omega_i$.  
 Consider  $\Omega_j=(\mu^i_j)^{-1}(\Omega_i)$.  Analyzing Procedures I-IV we deduce that $\dim \Omega_j=2$, 
  for Procedures I and IV we have that   $\Omega_{j+1}=\Omega_j$, and  for Procedures II and III 
  we have that   $\Omega_{j+1}$ is obtained 
 from $\Omega_j$ by taking a triangulation of $\Omega_j$ and replacing each $2$-simplex $\Delta$ of
 $\Omega_j$ by 
 by a mapping cylinder of a map $\partial\Delta \lo S^1$  of degree $k_j$ attached to $\partial \Delta$.
 Note that $k_j \geq k_{i+1}$ for $j >i$. Also  note that $\mu^i_j$ is $1$-to-$1$
 over  the $1$-skeleton $M^{(1)}_i$ of $M_i$ and hence
  $M_i^{(1)}$  can be considered as  a subset 
  $M^{(1)}_i \subset M^{(1)}_j$  of the $1$-skeleton of $M_j$
  and  then  $\partial \Omega_i $ considered as a subset of $\Omega_j$  is homologous
  to $0$ in $H_1(\Omega_j; \z_{p^{k_{i+1}}})$, $j>i$. 
  
  Consider the trivialization 
  $g^{-1}_j(\Omega_j)=\Omega_j \times S^1$ of $L_j$ over  $\Omega_j$
   induced by the trivialization of 
  $g_i^{-1}(\Omega_i)=\Omega_i \times S^1$ 
   and consider the section  $s_j : M^{(1)}_j \lo L_j$ which is $s^*_j$ followed by $G^*_j$.
   Note that for every $2$-simplex $\Delta$ of $\Omega_j$, the map $s_j$ restricted to  $\partial \Delta$ and
   followed by the projection of  $g^{-1}_j(\Omega_j)=\Omega_j \times S^1$ to $S^1$ is of
   degree $p^{k_{i+1}}$. Then, since $\partial \Omega_i $  is homologous
  to $0$ in $H_1(\Omega_j; \z_{p^{k_{i+1}}})$, we get that $s_j$ restricted to  $\partial \Omega_i$ and
   followed by the projection of  $g^{-1}_j(\Omega_j)=\Omega_j \times S^1$ to $S^1$ is of
   degree $p^{k_{i+1}}$ as well. Thus we get that the section
   $s_i=\mu^i_j \circ s_j |  M^{(1)}_i  : M^{(1)}_i \lo L_i$ is of degree $p^{k_{i+1}}$ and this violates  our choice of 
   $k_{i+1}$ in Procedure I. 
 $\black$

\end{section}

\begin{section}{Proof of Theorem \ref{dim-z=2}}
Theorem \ref{dim-z=2} follows from Theorem \ref{observation} and the following proposition.
\begin{proposition}
\label{proposition-dim-z=2}
There is an infinite dimensional compactum $Y$ with $\p Y=1$ and $\dim_\z Y=2$ such that
$H^2(Y; \z)$ does not contain a subgroup isomorphic to $\z_{p^\infty}$. 
\end{proposition}

Let $K$ be a CW-complex and $N$ a connected CW-complex. By $ \map (K,N)$ we denote
the space of  pointed maps from $K$ to $N$ with the compact-open topology. 
We will write  $\map (K,n)\cong 0$ if $\map (K,N)$ is weakly homotopy equivalent to a point.
Note that $\map (K,n)\cong 0$ if and only if for every $n\geq 0$,  every  map from $\Sigma^n K$ to  $N$
is null-homotopic.

In the proof of Proposition \ref{proposition-dim-z=2} we will use the following facts.

\begin{theorem}{\rm (Miller’s theorem (The Sullivan conjecture)) \cite{sullivan}}.
\label{sullivan}
Let $G$ be a finite group and $N$ a connected  finite CW-complex. 
Then $\map (K(G,1), N)\cong 0$.
\end{theorem}

\begin{proposition}{\rm \cite{levin-canada}}
\label{levin-canada}
 Let $M$ be a countable CW-complex,  
$N$  a connected CW-complex whose
homotopy groups are finite, $A_N$  a subcomplex of $M$,   and   $ f : A_N  \lo  N$ a map that cannot be
continuously extended over $M$  then there exists a finite subcomplex $M_N$  of $M $ such that 
$f_N | A_N\cap M_N : A_N\cap  M_N \lo N$ 
cannot be continuously extended over $ M_N$.
\end{proposition}

\begin{proposition}
\label{sullivan-1/p}
 Let  $N$ be a  CW-complex whose homotpy groups are finite $p$-torsion groups.
 Then
$\map( K(\z[\frac{1}{p}],1), N)\cong 0$. 
\end{proposition}
{\bf Proof.} Note that $\Sigma^n K(\z[\frac{1}{p}], 1)$ can be represented as the infinite telescope of
a map $S^n \lo S^n$ of degree $p$. Then $\pi_n (\Sigma^n K(\z[\frac{1}{p}], 1))=\z[\frac{1}{p}]$
and every  finite subtelesope of $K(\z[\frac{1}{p}], 1)$ 
 is homotopy equivalent to $S^n$. Since $\pi_n(N)$ is a $p$-torsion group we have that
  every map  $  f: \Sigma^n K(\z[\frac{1}{p}], 1))\lo N$ sends 
 $\pi_n (\Sigma^n K(\z[\frac{1}{p}], 1)) $ to $0$ in $\pi_n(N)$ and hence $f$ is null-homotopic on
 every finite subtelescope of $\Sigma^n K(\z[\frac{1}{p}], 1)$. Then
 $f$ extends over $\Sigma(\Sigma^n K(\z[\frac{1}{p}], 1))$ since otherwise  it would contradict
   Proposition \ref{levin-canada}. Thus 
 $f$ is null-homotopic.
 $\black$

\begin{proposition}
\label{partial-maps-more}
Consider Construction \ref{partial-maps}.  Let $N$ be  a connected CW-complex such that
  $\map(K,N)\cong 0$ and 
 $f_N : A_N \lo N$ is a map from a subcomplex $A_N$ of $M$  such that $f_N$ does not extend over $M$.
Then $\mu$ restricted to $\mu^{-1}(A_N)$ and followed by $f_N$ does not extend over $M'$.
\end{proposition}
{\bf Proof.} Let  $M*$ by the space obtained from $M'$ by collapsing  the fibers of $\mu$ over $N$ to singletons
and  $\mu^* : M^* \lo M$ the induced map.  Thus we can consider $N$ as a subcomplex of $M^*$ and 
aiming at a contradiction assume that  $f_N$  extends
to $f^* : M^*  \lo N$. 

Consider the first barycentric subdivision $\beta M$ of the triangulation of $M$.
Note that for ever $0$-dimensional simplex (vertex) $\Delta$ of $\beta M$,
$(\mu^*)^{-1}(\Delta)$  is either a singleton or homeomorphic to $K$. Then $f^*$ is null-homotopic
on $(\mu^*)^{-1}(\Delta)$. Thus we can replace $M^*$ by the space  $M^*_0$  obtained from $M^*$ by collapsing
$(\mu^*)^{-1}(\Delta)$  to  singletons for every
$0$-simplex $\Delta$ of $\beta M$  and $\mu^*$ and $f^*$  by the induced maps $\mu^*_0$ and $f^*_0$
to $M$ and
$N$ respectively   and assume that
$\mu_0^*$ is $1$-to-$1$ over the $0$-simplexes of $\beta M$.

Now consider a $1$-simplex $\Delta$ of $\beta M$. Then $(\mu^*_0)^{-1}(\Delta)$ is either contractible 
or homotopy equivalent to $\Sigma K$ and hence  $f^*_0$ is null-homotopic on
$(\mu^*_0)^{-1}(\Delta)$. Thus $f^*_0$ factors up to homotopy  through the space obtained from
$M^*_0$ by collapsing the fibers of $\mu^*_0$ over the simplex $\Delta$. Doing that consecutively 
for all the  $1$-simplexes of $\beta M$  we obtain the space $M^*_1$ and  the maps $\mu^*_1 : M^*_1 \lo M$
and $f^*_1: M^*_1 \lo N$  induced by $\mu^*_0$ and $f^*_0$ respectively such that 
$\mu^*_1$ is $1$-to-$1$ over the $1$-simplexes of $\beta M$.

Procced by induction and finally get for $m=\dim M$ that $M^*_m=M$ and
$f^*_m$ extends $f_N$ that  contradicts the assumptions of the proposition. $\black$
\\
\\
{\bf Proof of Proposition \ref{proposition-dim-z=2}.}
We will construct   $Y$ as the inverse limit of $(n+2)$-dimensional  finite  simplicial complexes $M_i$ and
 combinatorial
 bonding maps $\mu_{i+1}: M_{i+1} \lo M_i$. In order to show that 
 $Y$  has the required properties
 we consider
 for each  $i$   a subcomplex $A_i$ of
 $M_i$ such that  $A_{i+1}=\mu^{-1}_{i+1}(A_i)$ and 
  a  natural numbers $k_i$ 
 such that $k_{i+1} \geq k_i$.
Set  $A_0$ to be  a  Moore space $M(\z_p, 2)$,   $M_0$ the cone over $A_0$  
 and  $k_0=1$.
 We denote  $\mu^i_j =\mu_j \circ \dots \circ \mu_{i+1}: M_j \lo M_i$ with $\mu^i_i$ being the identity map
 and require that

\begin{enumerate}[start=1, label={ (\arabic*)}]

\item $\mu^0_i | A_i  :  A_i \lo A_0=M(\z_p, 2)$ does not  extend over $M_i$ as a map to $M(\z_p, 2)$.
  
\end{enumerate}
Clearly   (1)   holds for $i=0$.
Assuming that  the construction is completed for $i$  and we proceed to $i+1$ performing one the following 
procedures. \\

 {\bf Procedure I (taking care of non-flexible  bundles).} Let $g : L \lo M_i$ be a circle bundle which is not
   $p$-flexible. Then take  any  natural number $k$ such that  there is no section over $M_i^{(1)}$ of degree $p^k$.
   Set $M_{i+1}=M_i$, $\mu_{i+1}$=the identity map and  $k_{i+1}= \max\{k, k_i\}$. 
   Clearly (1) holds for $i+1$.\\
   
   {\bf Procedure II (taking care of flexible bundles).}
   Note that the homotopy groups of $M(\z_p, 2)$ are finite  $p$-torsion groups, denote by $t$  any  natural number such that 
  $t \geq k_i$ and  $p^t\pi_j(M(\z_p, 2))=0$ for every $j \leq m=\dim M_i$.
   Let $g : L \lo M_i$ be a flexible bundle. Apply \ref{flexible-2} with 
   $M=M_i, N= M(\z_p, 2), A_N=A_i, f_N =\mu^0_i|A_i : A_N =A_i \lo N=A_0=M(\z_p, 2)$ and
    $t$ and $m$ as above to
   construct $M_{i+1}=M' $ and $\mu_{i+1}=\mu'$, and   set $k_{i+1}=k_i$.
   By Proposition \ref{flexible-2-more}, we get that  (1) holds for $i+1$.\\
   
{\bf Procedure III (extending partial maps to a circle).} 
Let $A$ be a subcomplex of $M_i$ and $f : A \lo K(\z [\frac{1}{p}], 1)$ a map.
Recall that $K(\z [\frac{1}{p}], 1)$ can be represented as the infinite telescope of the $\z_p$-covering map 
$S^1 \lo S^1$. 
Extend $f$ over  the $1$-skeleton $M^{(1)}$ of $M$  and replace $f$ by  a cellular map
homotopic  to $f$ such  that $f(A)=S^1\subset  K(\z [\frac{1}{p}], 1)$  and $f$ as a map to  $S^1$ is 
of degree $k_i$.  Apply \ref{partial-maps} with $M=M_i$ and
   $K=K(\z [\frac{1}{p}], 1)$  
  to  construct $M_{i+1}=M' $ and  $\mu_{i+1}=\mu$. 
  By Propositions \ref{sullivan-1/p} and  \ref{flexible-2-more} we get  that (1) holds for $i+1$.
  By Proposition \ref{levin-canada}  we can replace 
  in Construction \ref{partial-maps} the complex $K= K(\z [\frac{1}{p}],1)$ by 
  a finite subcomplex of $K(\z [\frac{1}{p}],1)$
  containing $f(A)$ and still preserve (1). Thus we
    assume that $M'$ is a finite CW-complex and set
   $k_{i+1}=k_i$.\\
   
  {\bf  Procedure  IV (extending partial maps to $K(\z_{p^\infty}, 1)$).} 
  Let $A$ be a subcomplex of $M_i$
  and $f: A \lo K(\z_{p^\infty}, 1)$ a   map.  Extend $f$ over the $1$-skeleton of $M_i$.
  Represent $K(\z_{p^\infty}, 1)$ as the infinite telescope of the maps 
  $K(\z_{p^n}, 1)\lo K(\z_{p^{n+1}}, 1)$ induced by the monomorphisms $\z_{p^n} \lo \z_{p^{n+1}}$
  and replace $f$ by a homotopic  cellular map   $f: A \lo K(\z_{p^n}, 1) \subset K(\z_{p^\infty}, 1)$ 
  such that $f(M_i^{(1)})=S^1$ and $f$ restricted to $M^{(1)}_i$ and considered as a map to $S^1$
  is of degree $p^{k_i}$.
  Apply \ref{partial-maps} with $M=M_i$ and $K=K(\z_{p^n}, 1)$ to 
   construct $M_{i+1}=M' $ and $\mu_{i+1}=\mu$. 
    By Theorem \ref{sullivan}  and Proposition \ref{partial-maps-more} we get  that (1) holds for $i+1$.
     By Proposition \ref{levin-canada}  we can replace 
  in Construction \ref{partial-maps} the complex $K= K(\z_{p^n}, 1)$ by 
  a finite subcomplex of $K(\z_{p^n}, 1)$
  containing $f(A)$ and still preserve (1). Thus we
    assume that $M'$ is a finite CW-complex and set
   $k_{i+1}=k_i$.
   \\
   \\

Let   $Y=\invlim M_i$. Clearly (2) implies that  $\dim Y >2$.

 Assume that $Y$ is constructed as described in \ref{inverse-limit}. 
 The objects that we pick out for each $M_i$ are
 all the  (non-isomorphic) circle bundles  over $M_i$ and 
 countably many maps from subcomplexes of $M_i$ to $K(\z[\frac{1}{p}], 1) $ and $K(\z_{p^\infty}, 1)$ representing
 all possible maps up to homotopy.  The procedures that we perform are Procedures I-IV.
 
 Assuming that on each $M_i$ we fix a sufficiently finite triangulation, we  get that:
 
 $\bullet$  Procedure III implies 
  that \ext$ \leq K(\z[\frac{1}{p}], 1)$ and hence $\p Y \leq 1$ and
 
 $\bullet$ Procedure IV implies  that
 \ext$ Y \leq K(\z_{p^\infty}, 1)$ and hence
  $\dim_{\z_p} Y \leq n+1$.
 
 Then $\dim Y>2$ and Bockstein inequalities
 imply that $\dim_\z  Y =2$ and hence $Y$ is infinite dimensional.
 
  Let us show that $H^2(Y; \z)$ does not contain a subgroup isomorphic to $\z_{p^\infty}$. Since
  $\p Y =1$ we have that  $H^2(Y; \z)\otimes \z[\frac{1}{p}]=0$ and hence   $H^2(Y; \z)$ is a $p$-torsion group.
  Thus we need to show  $H^2(Y; \z)$ does not contain a non-trivial
   element $\alpha$ that is divisible by $p^k$ for every $k$.
  Aiming at a contradiction assume that such $\alpha$ does exist. Consider the circle bundle
  $g_Y : Z \lo Y$ corresponding to $\alpha$. Then Construction \ref{inverse-limit}
  guarantees that at a certain step $i$
of the construction we will consider a circle bundle $g : L\lo M_i$ such that
$g_Y$ is the pull-back of $g$ under the projection $\mu^i: Y \lo Y_i$.

 If $g$ is flexible then we apply Procedure II and,
by Construction \ref{flexible-1}, 
the pull-back bundle $g_{i+1} : L_{i+1} \lo M_{i+1}$ of $g$ under the map $\mu_{i+1}$ admits
 a section and hence both $g_{i+1}$  and  $g_Y$ are trivial. Thus we arrive at a  contradiction with $\alpha\neq 0$.
 
 If $g: L \lo M_i$ is non flexible then we apply Procedure I. Denote by $g_j : L_j \lo M_j, j\geq i$, 
 the pull-back bundle of $g$   to $M_j$ under 
 $\mu^i_j $
 and   $g_i=g$ and $L_i=L$. Denote by $G^i_j : L_j \lo L_i$ the induced map.
 Recall that  $\alpha$ is divisible by $p^k$  for every  $k$. Then for a sufficiently large $j$
  there is a circle bundle $g^*_j : L^*_j \lo M_j$ such that $g_j : L_j \lo M_j$ is induced  by 
  $g^*_j$ and $\z_{p^{k_{i+1}}}$,  see Lemma \ref{circle}.  Let 
  $G_j^* : L^*_j \lo L_j$ be a bundle map induced by the homomorphism $S^1 \lo S^1 / \z_{p^{k_{i+1}}}$,  
Take a section $s^*_j : M^{(1)}_j \lo L^*_j$.

 Let $\Omega$  be   a $2$-dimensional subcomplex of   $M_j$. Denote by 
 $\mu^{-\#}_{j+1}(\Omega)$ the subcomplex of  $\mu^{-1}_{j+1}(\Omega)$ which is the closure of 
 the following set
 
  \ \  \ $\mu^{-1}_{j+1}(\Omega)\setminus \{$
  the fibers of $\mu_{j+1}$ over the barycenters of the $2$-simplexes of  $\Omega \}$.
  \\
 Consider a $2$-simplex $\Omega_i$ of $M_i$
   and  let  $g^{-1}_i(\Omega_i)=\Omega_i \times S^1$ be
  a trivialization of $L_i$ over $\Omega_i$.   For every $j>i$  define  $\Omega_j$ 
    by  $\Omega_{i+1}=\mu^{-\#}_{i}(\Omega_i),  \Omega_{i+2}=
    \mu^{-\#}_{i+2}(\Omega_{i+1}), \dots, 
  \Omega_{j}=\mu^{-\#}_{j}(\Omega_{j-1})$.

  Analyzing Procedures I-IV we deduce that $\dim \Omega_j=2$, 
  for Procedure I  we have that   $\Omega_{j+1}=\Omega_j$, and  for Procedures II, III and IV
  we have that   $\Omega_{j+1}$ is obtained 
 from $\Omega_j$ by taking a triangulation of $\Omega_j$ and replacing each $2$-simplex $\Delta$ of
 $\Omega_j$
 by a mapping cylinder of a map $\partial\Delta \lo S^1$  of degree $k_t$ attached to $\partial \Delta$.
 Note that $k_j \geq k_{i+1}$ for $j >i$. Also  note that $\mu^i_j$ is $1$-to-$1$
 over  the $1$-skeleton $M^{(1)}_i$ of $M_i$ and hence
  $M_i^{(1)}$  can be considered as  a subset 
  $M^{(1)}_i \subset M^{(1)}_j$  of the $1$-skeleton of $M_j$
  and  then  $\partial \Omega_i $ considered as a subset of $\Omega_j$  is homologous
  to $0$ in $H_1(\Omega_j; \z_{p^{k_{i+1}}})$, $j>i$. 
  
  Consider the trivialization 
  $g^{-1}_j(\Omega_j)=\Omega_j \times S^1$ of $L_j$ over  $\Omega_j$
   induced by the trivialization of 
  $g_i^{-1}(\Omega_i)=\Omega_i \times S^1$ 
   and consider the section  $s_j : M^{(1)}_j \lo L_j$ which is $s^*_j$ followed by $G^*_j$.
   Note that for every $2$-simplex $\Delta$ of $\Omega_j$, the map $s_j$ restricted to  $\partial \Delta$ and
   followed by the projection of  $g^{-1}_j(\Omega_j)=\Omega_j \times S^1$ to $S^1$ is of
   degree $p^{k_{i+1}}$. Then, since $\partial \Omega_i $  is homologous
  to $0$ in $H_1(\Omega_j; \z_{p^{k_{i+1}}})$, we get that $s_j$ restricted to  $\partial \Omega_i$ and
   followed by the projection of  $g^{-1}_j(\Omega_j)=\Omega_j \times S^1$ to $S^1$ is of
   degree $p^{k_{i+1}}$ as well. Thus we get that the section
   $s_i=\mu^i_j \circ s_j |  M^{(1)}_i  : M^{(1)}_i \lo L_i$ is of degree $p^{k_{i+1}}$ and this violates  our choice of 
   $k_{i+1}$ in Procedure I. 
 $\black$

\end{section}

Michael Levin\\
Department of Mathematics\\
Ben Gurion University of the Negev\\
P.O.B. 653\\
Be'er Sheva 84105, ISRAEL  \\
 mlevine@math.bgu.ac.il\\\\
\end{document}